\documentclass[preprint,review,sort&compress,3p]{elsarticle}
\usepackage{xcolor}
\usepackage{hyperref}
\hypersetup{
    colorlinks=false,
    linkcolor=blue,
    urlcolor=cyan,
}
\usepackage{colortbl} % 
\usepackage{float}
\usepackage[utf8]{inputenc} % allow utf-8 input
\usepackage{threeparttable}
\usepackage{diagbox}
\usepackage{mathrsfs}
\usepackage{amsfonts}
\usepackage{cancel}%delete or cancel
\usepackage{graphicx}
\usepackage{amssymb}
\usepackage{hyperref}
\usepackage{amsthm}
\usepackage{amsmath}
\usepackage{epstopdf}
\usepackage{verbatim}
\usepackage{booktabs}%
\usepackage{flafter}
\usepackage{algorithm}
\usepackage{algorithmicx}
\usepackage{algpseudocode}
\usepackage{makecell}
\usepackage{caption,subcaption}
\usepackage{multirow}
\usepackage{appendix}
\usepackage{float}
\usepackage{array, caption, threeparttable}
\usepackage[font=small,labelfont=bf,labelsep=none]{caption}
\numberwithin{equation}{section}

\journal{ }

\begin{document}
\begin{frontmatter}
\title{U-WNO: U-Net Enhanced Wavelet Neural Operator for Solving Parametric Partial Differential Equations}
%\thanks[b]{\small Supported by National Natural Science Foundation of China (11101071, 11271001, 51175443) and Fundamental Research Funds for the Central Universities (ZYGX2016J131, ZYGX2016J138).}
\cortext[cor1]{Corresponding author}
\author[a]{Wei-Min Lei}
\author[a]{Hou-Biao Li \corref{cor1}}
\ead{lihoubiao0189@163.com}

\address[a]{School of Mathematical Sciences, University of Electronic Science and Technology of China, Chengdu, 611731, P. R. China}
%\date{November 17, 2023}

\begin{abstract}
High-frequency features are critical in multiscale phenomena such as turbulent flows and phase transitions, since they encode essential physical information. The recently proposed Wavelet Neural Operator (WNO) utilizes wavelets' time-frequency localization to capture spatial manifolds effectively. While its factorization strategy improves noise robustness, it suffers from high-frequency information loss caused by finite-scale wavelet decomposition. In this study, a new U-WNO network architecture is proposed. It incorporates the U-Net path and residual shortcut into the wavelet layer to enhance the extraction of high-frequency features and improve the learning of spatial manifolds. Furthermore, we introduce an adaptive activation mechanism to mitigate spectral bias through trainable slope parameters. Extensive benchmarks across seven PDE families (Burgers, Darcy flow, Navier-Stokes, etc.) show that U-WNO achieves 45--83\% error reduction compared to baseline WNO, with mean $L_2$ relative errors ranging from 0.043\% to 1.56\%. This architecture establishes a framework combining multiresolution analysis with deep feature learning, addressing the spectral-spatial tradeoff in operator learning. Code and data used are available on \url{https://github.com/WeiminLei/U-WNO.git}.
\end{abstract}

\end{frontmatter}

\section{Introduction}

In recent years, Neural Networks (NN) \cite{minar2018recent} have been widely used in computational science to solve classical applied mathematics problems. This is due to their strong approximation capabilities and ability to handle nonlinear systems. Natural systems are governed by conservation and constitutive laws. Modeling complex constitutive laws using PDEs is a popular approach in various fields of science and engineering disciplines \cite{debnath2005nonlinear,evans2022partial,jones2009differential}. In most cases, the preferred methods for solving these PDEs, due to the lack of analytical solutions, are the Finite Element Method (FEM) \cite{taylor1998finite}, isogeometric analysis \cite{cottrell2009isogeometric}, the Finite Difference Method (FDM) \cite{zhang2009finite}, and the Finite Volume Method (FVM) \cite{eymard2000finite}. Sometimes, in science and engineering, it becomes necessary to solve PDEs with different parameters, such as domain geometries, source functions, and initial or boundary conditions (ICs and BCs). These types of PDEs are often referred to as parametric PDEs. However, the computational cost of the traditional numerical schemes described above is high in such cases, as they require independent calculations for each combination of input parameters. In recent years, several approaches using Machine Learning (ML) \cite{sah2020machine} have been proposed as faster alternatives to numerical simulations \cite{tahmasebi2020machine}. For example, Physics-Informed Neural Networks (PINNs) \cite{raissi2019physics}, Deep Galerkin methods (DGM) \cite{sirignano2018dgm}, and Deep Ritz method (DRM) \cite{yu2018deep} satisfy given governing equations and boundary conditions by constraining the output of NN. However, when solving parametric PDEs using these methods, it is impractical to train a model individually for each PDE parameter. Therefore, deep learning methods for solving parametric PDEs are often considered.

Deep learning (DL) methods for solving parametric partial differential equations are usually divided into two categories:

(1) By combining meta-learning \cite{vanschoren2018meta} etc., the aforementioned schemes are used to solve parametric PDEs. For example, some researchers have proposed using transfer learning \cite{tan2018survey} or meta-learning to fine-tune these schemes on the original trained network for rapid generalization to new parameters. Ref. \cite{yu2018deep} recommends a transfer learning approach that utilizes a model trained for one task as the initial model for training another task, which can accelerate the training of PINNs. This method is effective in the case of small parameter variation ranges, but not in the case of large parameter variation ranges. Therefore, based on the minimum energy path theory, the PDE parameters are set as the trainable parameters of the NN, and the low-loss path is "preset" to ensure training stability and guide the migration of the remaining network weights \cite{goswami2020transfer}. Based on the idea of Meta-Learning, a Meta-Learning PINN (Meta-PINN) \cite{zhong2023accelerating} is designed for plasma simulation acceleration. Some researchers have combined the Reptile algorithm \cite{nichol2018reptile} with PINNs to accelerate the solution of parametric PDEs \cite{liu2022novel}. Meta-MgNet \cite{chen2022meta} is a work that treats solving parametric PDEs as a Meta-Learning problem. Meta-MgNet leverages the similarity between tasks to adaptively generate a good smoothing operator, which speeds up the solution process. Inspired by the Reduced Basis Method (RBM) \cite{bai2005reduced}, Dr. Yanlai Chen designed a pre-trained PINN (GPT-PINN) for parametric PDEs \cite{chen2024gpt}.

(2) Another approach utilizes a neural network (NN) to learn how to map the parameters of partial differential equations (PDEs) to their solutions. For example, the Deep Operator Network (DeepONet) \cite{lu2019deeponet} introduces two nonlinear operator theories that can be implemented using an NN architecture. DeepONet uses a trunk net and a branch net to encode the parameters and spatial positions of the PDEs, respectively. It combines the outputs of these two networks to obtain the solution of the PDEs. In addition, the Fourier Neural Operator (FNO) \cite{li2020fourier} directly parameterizes the integral kernel in the Fourier space, resulting in an NN architecture with stronger representation capabilities. By approximating the Green function using the powerful fitting ability of the NN, it transforms the integration over the entire control domain of the PDEs into the integration over the node domain. Recently, a graph neural operator (GNO) for parametric PDEs is proposed \cite{li2020neural}. To reduce the complexity of the calculation, the Multipolar Graph Neural Operator (MGNO) \cite{li2020multipole} is introduced, which reduces the approximation cost of the kernel function. However, FNO suffers from a lack of spatial resolution due to the frequency localization of the FFT basis function \cite{bracewell1989fourier}, leading to underperformance in PDEs with complex geometries. To address this limitation, the Wavelet Neural Operator (WNO) \cite{tripura2023wavelet,tripura2022wavelet} is proposed, which uses the wavelet transform \cite{zhang2017time} to learn the spatial frequency positioning and mode changes in the input. WNO can be seen as a more general case of the previously proposed FNO. While WNO excels at learning complex nonlinear differential operators regardless of geometric structure, it restricts the parameterization of neural network weights to high-level wavelet decomposition. This helps mitigate the effects of data noise but hinders the extraction of high-frequency features, leading to decreased accuracy. To address this issue, we propose an enhanced version of WNO called U-WNO:
\begin{itemize}
\item Firstly, the U-Net \cite{ronneberger2015u} path is incorporated into the wavelet layer of the WNO, which samples the data at various scales. This enables U-WNO to reconstruct information at different frequency scales. Based on the U-Net path, the U-WNO can extract and fuse potential features more effectively.
\item By introducing residual shortcuts into the wavelet layer, we ensure hierarchical propagation of complete input function information across scales. This enhances the performance of neural operators in spatial domains.  
\item Finally, because of the Spectral Bias \cite{rahaman2019spectral} of the NN, we apply the Adaptive Activation Function \cite{jagtap2020adaptive, jagtap2022deep, jagtap2020locally}  to improve the convergence speed and accuracy of the network.
\end{itemize}
\par
 
The structure of this article is as follows. In Section 2, we introduced the general content of neural operators. We describe our proposed methodology in detail in Section 3. In Section 4, we conduct numerical simulations to demonstrate the effectiveness of our method. Finally, We conclude with some thoughts on possible directions for further research in Section 5.

\section{Problem setting and neural operators}
An operator is a mapping between an infinitely dimensional input and output function space. Let $D\in\mathbb{R}^{d}$ be the d-dimensional domain of the boundary $\partial D$. For the fixed domains $D=\left( a,b \right)$ and $x\in D$, consider a PDE that maps the input function space, i.e., the space containing the source term $f(x,t):D\mapsto\mathbb{R}$, the initial condition $u(x,0):D\mapsto\mathbb{R}$, and the boundary condition $u(\partial D,t):D\mapsto\mathbb{R}$ to the solution space $u(x,t):D\mapsto\mathbb{R}$, where $t$ is the time coordinate. We aim to learn an operator from the input function to the solution space $u\left( x,t \right)$. Now, we define two complete norm vector spaces $(\mathcal{A},\|.\|_{\mathcal{A}})$ and $(\mathcal{U},\|.\|_{\mathcal{U}})$ that take the values of a given partial derivative number in $\mathbb{R}^{d_{a}}$ and $\mathbb{R}^{d_{u}}$. These function spaces, often referred to as Banach spaces, are represented as $\mathcal{A}:=\mathcal{C}(D;\mathbb{R}^{d_{w}})$ and $\mathcal{A}:=\mathcal{C}(D;\mathbb{R}^{d_{w}})$. For the given domains $D\in\mathbb{R}^{d}$ and the canonical spaces $(\mathcal{A},\|.\|_{\mathcal{A}})$ and $(\mathcal{U},\|.\|_{\mathcal{U}})$, our goal is to learn the nonlinear operators $\mathcal{D}:\mathcal{A}\mapsto\mathcal{U}$. Once successful, a series of PDE solutions can be obtained for different input parameter sets $a\in \mathcal{A}$. If there exists a differential operator $\mathcal{N}$ that maps the function spaces to null space, i.e., $\mathcal{N}: \mathcal{A} \times \mapsto \mathcal{O}$, where $\mathcal{O}$ is the null space, then the parametric PDEs take the form:
\begin{align}\label{eq:2.1}
	\mathcal{N}\left( a,u \right) =0,\ \text{ in }D\subset \mathbb{R}^d.
\end{align}
The PDE is defined on a $d$-dimensional bounded domain, $D \in \mathbb{R}^{d}$, with the boundary $\partial D$, where the boundary condition is expressed as:
\begin{align}\label{eq:2.2}
	u = g, \,\, \text { in } \partial D .
\end{align}
The input and output functions in the function spaces $\mathcal{A}$ and $\mathcal{U}$ are represented by $a\ :\ D\mapsto a\left( x\in D \right) \in \mathbb{R}^{d_a}$ and $u\ :\ D\mapsto u\left( x\in D \right) \in \mathbb{R}^{d_u}$. Considering $u_{j}=\mathcal{D}(a_{j})$, we have access to $N$ pairs of $\{(a_j,u_j)\}_{j=1}^N$. Then, our goal is to approximate $D$ by a neural network (NN) as
\begin{align}\label{eq:2.3} \mathcal{D}:\mathcal{A}\times\boldsymbol{\theta}_{NN}\mapsto\mathcal{U},\quad\boldsymbol{\theta}_{NN}=\left\{\mathcal{W}_{NN},\boldsymbol{b}_{NN}\right\},
\end{align}
where $\boldsymbol{\theta }_{NN}$ denotes the finite-dimensional parameter domain of the NN. Although the input and output functions $a_j$ and $u_j$ are continuous, for numerical execution, we discretize the field $D$ with $n_D$ points $\{x_{p}\}_{p=1}^{n_{D}}$. Thus, for $N$ sets of input and output pairs, supposing we have access to the dataset $\{a_{j}\in\mathbb{R}^{n_{D}\times d_{a}}$, $ u_{j}\in\mathbb{R}^{n_{D}\times d_{u}}\}_{i=1}^{N}$. The goal is to develop a network capable of utilizing data to learn the operator $\mathcal{D}$. For featuring a multidimensional kernel convolution, the input $a\left( x \right) \in \mathcal{A}$ is lifted to a higher-dimensional space through the local transformation $P(a(x)):\mathbb{R}^{d_a}\mapsto\mathbb{R}^{d_v}$ and represent the higher-dimensional space as $v_{0}(x)=P(a(x))$, where $v_{0}(x)$ takes the value in $\mathbb{R}^{d_v}$. In addition, $L$ steps are required to converge. We apply $L$ iterations $v_{i+1}=G(v_{i})$ for $i=0,\ldots,L$ on $v_0\left( x \right)$, where the function $G:\mathbb{R}^{d_{v}}\mapsto\mathcal{R}^{d_{v}}$ takes its value in $\mathbb{R}^{d_{v}}$. Once the iteration is done, we define another local transformation $Q(v_{L}(x)):\mathbb{R}^{d_{v}} \mapsto \mathbb{R}^{d_{u}}$ and convert the output $v_L\left( x \right)$ back to the solution space $u(x)=Q\left(v_{L}(x)\right)$. $v_{j+1}=G(v_{j})$ is defined as
\begin{align}\label{eq:2.4}
	v_{j+1}(x):=\sigma\left(\left(K(a;\phi)*v_{j}\right)(x)+Wv_{j}(x)\right),\quad x\in D,\quad j\in[0,L],
\end{align}
where $\sigma(\cdot):\mathbb{R}\to\mathbb{R}$ is a nonlinear activation function, $W:\mathbb{R}^{d_{v}}\to\mathbb{R}^{d_{v}}$ is a linear transformation, and $K:\mathcal{A}\times \boldsymbol{\theta} \rightarrow \mathcal{L}\left( \mathcal{U},\mathcal{U} \right)$ is an integral operator on $\mathcal{C}(D;\mathbb{R}^{d_{v}})$. In the NN setup, we set $K(a;\phi)$ to represent a kernel integral operator parameterized by $\phi \in \boldsymbol{\theta }$. The kernel integral operator $K(a;\phi)$ is defined as,
\begin{align}\label{eq:2.5} \left(K(a;\phi)*v_{j}\right)(x):=\int_{D\in\mathbb{R}^{d}}k\left(a(x),x,y;\phi\right)v_{j}(y)\mathrm{d}y;\quad x\in D,\quad j\in[0,L],
\end{align}
where $k_\phi: \mathbb{R}^{2d+d_a} \mapsto \mathbb{R}^{d_v\times d_v}$ is an NN parameterized by $\begin{matrix}{\phi}&{\in}&{\boldsymbol{\theta}}\\\end{matrix}$. Suppose some kernel $k \in \mathcal{C}\left(D;\mathbb{R}^{d_{v}}\right)$ then we learn $\kappa_{\phi}$ from the data with the concept of degenerate kernels, so that,
\begin{align}\label{eq:2.6}
\|\begin{pmatrix}k_\phi(a)-k(a)\end{pmatrix}v\|\leq\|\begin{pmatrix}k_\phi(a)-k(a)\end{pmatrix}\|\|v\|,
\end{align}
and $k_{\phi}\rightarrow k$. Since the kernel $k\in\mathcal{C}(D)$ and $K \in \mathcal{L}(\mathcal{C}(D;\mathbb{R}^{d_{v}}))$ together define the infinite-dimensional space, Eq.(\ref{eq:2.4}) can learn the mapping of any infinite-dimensional function space. When it comes to approximating nonlinear operators, the structure of Eq.(\ref{eq:2.4}) follows that of a traditional neural network. To train the network, we define an appropriate loss function, $\mathcal{L}=\mathcal{U}\times\mathcal{U}$ as follows:
\begin{align}\label{eq:2.7}
\theta_{NN}=\arg\min_{\theta_{NN}}\mathcal{L}\left(\mathcal{D}(a),\mathcal{D}(a,\theta)\right).
\end{align}
Then, the Eq. (\ref{eq:2.7}) is minimized in a data-driven way to obtain the best parameters of the network. Finally, we obtain the neural operator.

\section{Proposed Methodology}\label{sec3}
\vspace{-2mm}
\subsection{Wavelet Neural Operator $^{\cite{tripura2023wavelet}}$}\label{sec3.1}

To delete the term $a\left(x,y \right)$ from Eq. (\ref{eq:2.5}), change the kernel to $k_{\phi}(x-y)$. We rewrite Eq.(\ref{eq:2.5}) as follows:
\begin{align}\label{eq:3.1}
\left(K(\phi)*v_{j}\right)(x):=\int_{D\in\mathbb{R}^{d}}k\left(x-y;\phi\right)v_{j}(y)dy,\quad x\in D,\quad j\in[0,L].
\end{align}
This is consistent with the convolution operator. The kernel $k_{\phi}$ is learned by parameterization in the wavelet domain. The continuous wavelet transform (CWT) is computationally expensive, so the discrete wavelet transform (DWT) is chosen by learning the kernel in the wavelet domain. The above convolution is then performed on the wavelet decomposition coefficients. Let $\psi(x)\in\mathcal{L}^2(\mathbb{R})$ be an orthogonal parent wavelet. Assuming that $\mathcal{W}(\Gamma)$ and $\mathcal{W}^{-1}(\Gamma_{W})$ are the forward and reverse discrete wavelet transforms of any function $\varGamma: D\mapsto \mathbb{R}^{d_v}$, the forward and reverse wavelet transforms of the function $\varGamma$ with scaling and shift parameters $m\in\mathbb{R}^{+}$ and $\tau\in\mathbb{R}$ are given by the following integral pairs,
\begin{align}\label{eq:3.2}
	\begin{aligned}
		&(\mathcal{W}\Gamma)_{j}(m,\tau) =\frac{1}{\sqrt{2^{m}}}\int_{D}\Gamma(x)\psi\left(\frac{x-\tau2^m}{2^m}\right)dx, \\
		&(\mathcal{W}^{-1}\Gamma_{w})_{j}(x) =\frac{1}{C_{\psi}\sqrt{2^{m}}}\int_{0}^{\infty}\int_{D}(\Gamma_{w})_{j}(m,\tau)\tilde{\psi}\left(\frac{x-\tau2^m}{2^m}\right)d\tau\frac{dm}{m^{2}}, 
	\end{aligned}
\end{align}
The term $C_{\psi}$ is called the allowable constant and has the range $0<C_{\psi}<\infty $. The expression of $C_{\psi}$ is
\begin{align}\label{eq:3.3}
C_{\psi}=2\pi\int_{D}\frac{\left|\psi(\omega)\right|^{2}}{\left|\omega\right|}d\omega.
\end{align}
Applying the convolution theorem, let $R_{\phi}=\mathcal{W}(k_{\phi})$, Eq.(\ref{eq:3.1}) on the wavelet domain is expressed as
\begin{align}\label{eq:3.4}
\left(\mathcal{K}(\phi)*v_{j}\right)(x)=\mathcal{W}^{-1}\left(R_{\phi}\cdot\mathcal{W}(v_{j})\right)(x),\quad x\in D.
\end{align}
$R_{\phi}$ is parameterized only at the highest decomposition level, so we have $\mathcal{W}\left(v_{t}\right)\in\mathbb{R}^{n/2^{m}\times d_{v}}$ in the m-level of DWT. The length of the wavelet coefficient is also affected by the number of vanishing moments of the orthogonal parent wavelet. Therefore, it is more appropriate to write $\mathcal{W}(v_{t},m)\in\mathbb{R}^{\zeta\times d_{v}}$, where the $\zeta$ wavelet coefficient supports the length of the width. Thus, direct parameterization of $R_{\phi}$ in the wavelet space is performed by integrating the parameter $\phi $ into the convolution of the value tensor of $(\zeta\times d_v\times d_v)$. To rewrite Eq.(\ref{eq:3.4}), we can write the convolution of the weight tensor $R\in\mathbb{R}^{\zeta\times d_v\times d_v}$ and $\mathcal{W}(v_{t},m)\in\mathbb{R}^{\zeta\times d_{v}}$ as,
\begin{align}\label{eq:3.5}
\left(R\cdot\mathcal{W}(v_j,m)\right)_{I_1,I_2}(x)=\sum_{I_3=1}^{d_v}R_{I_1,I_2,I_3}\mathcal{W}(v_j,m)_{I_1,I_3},\quad I_1\in[1,\zeta],\quad I_2,I_3\in d_v.
\end{align}
The WNO architecture consists of three steps:
\begin{itemize}
	\item The input function $a(x)$ is elevated to a high-dimensional space $v_0 (x) =P (a(x))$ by transforming $P$ with a fully connected NN (FNN).
	\item Application of wavelet layer: $v_{0}\mapsto v_{1}\mapsto...\mapsto v_{L}$,  where $v_{j}$ for $j = 0,1,...,L$ is a sequence of functions taking the value of the channel dimension $c$ in $\mathbb{R}^{c}$.
	\item Use the FNN transform $Q$ to project $v_{L}$ back to the original space $u(x)=Q(v_{L}(x))$. 
\end{itemize}
In each of the wavelet layers, there are
\begin{align}\label{eq:3.6}
	v_{j+1}(x):=\sigma\bigg(\big(\mathcal{K}v_{j}\big)(x)+W(v_{j}(x))\bigg),\forall x\in D,
\end{align}
where $\mathcal{K}$ is the kernel integral transform defined above and $W$ is a linear operator, both of which are learnable. $\sigma$ is an activation function that introduces nonlinearities to each wavelet layer.

\subsection{The Spectral Bias and Adaptive Activation Function}\label{sec3.2}
We first review the use of Neural Tangent Kernel (NTK) \cite{jacot2018neural} theory to study the Spectral Bias of NN \cite{wang2021eigenvector} and verify it empirically. Then the experiment confirms that the Spectral Bias of NN can be overcome by introducing the Adaptive Activation Function. Let $f(x,\boldsymbol{\theta})$ be a fully connected neural network (FNN) with weight $\boldsymbol{\theta}$ initialized by the Gaussian distribution $\mathcal{N}(0,1)$. Given a dataset $\{\boldsymbol{X}_{\mathrm{train}},\boldsymbol{Y}_{\mathrm{train}}\}$, where $\boldsymbol{X}_{\mathrm{train}}=(x_i)_{i=1}^N$ is the input and $\boldsymbol{Y}_{\mathrm{train}}=(y_i)_{i=1}^N$ is the output, we consider the mean square loss $\mathcal{L}(\theta)=\frac{1}{N}\sum_{i=1}^{N}|f(x_{i},\boldsymbol{\theta})-y_{i}|^{2}$ using a very small learning rate $\eta$. The neural tangent kernel operator $\boldsymbol{K}$ can be defined as follows:
\begin{align}\label{eq:3.7}
	K_{ij}=K(x_i,x_j)=\left\langle\frac{\partial f(x_i,\boldsymbol{\theta}))}{\partial\boldsymbol{\theta}},\frac{\partial f\left(x_j,\boldsymbol{\theta}\right))}{\partial\boldsymbol{\theta}}\right\rangle.
\end{align}

\begin{figure}[htbp]
	\centering
	\includegraphics[scale=0.25]{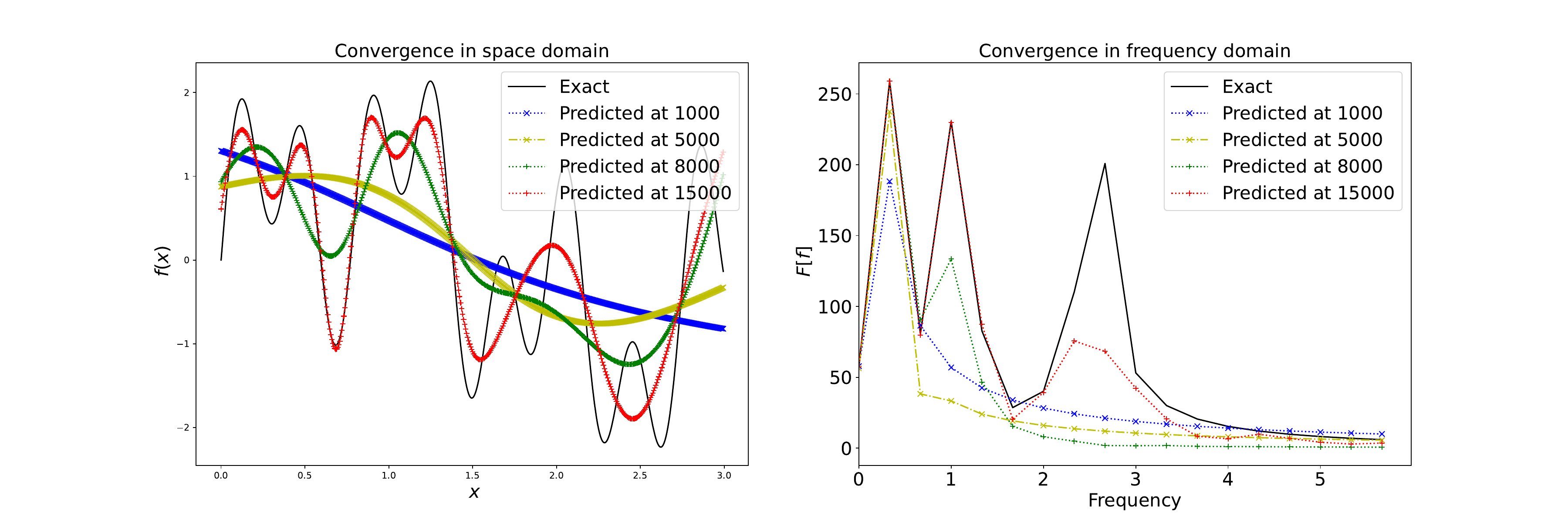}
	\caption{\; The convergence of FNN in the spatial domain (left) and frequency domain (right) when approaching $f(x)$.}
	\label{fig1}
\end{figure}

The NTK theory states that during gradient descent dynamics with infinitesimal learning rates, the kernel $K$ converges to a deterministic kernel $K^{*}$ and remains unchanged during training as the width of the network approaches infinity. Additionally, we have
\begin{align}\label{eq:3.8}
\frac{df(X_{\mathrm{train}},\theta(t))}{dt}\approx-K\cdot(f(\boldsymbol{X}_
{\mathrm{train}},\boldsymbol{\theta}(t))-\boldsymbol{Y}_{\mathrm{train}}),
\end{align}
$\boldsymbol{\theta}(t)$ denotes the parameters of the network in the iteration $t$, $f(\boldsymbol{X}_{\mathrm{train}},\boldsymbol{\theta}(t))=(f(x_i,\boldsymbol{\theta}(t))_{i=1}^N$. Then, it has the following form
\begin{align}\label{eq:3.9}
	f(\boldsymbol{X}_{\mathrm{train}},\boldsymbol{\theta}(t))\approx(I-e^{-\boldsymbol{K}t})\cdot\boldsymbol{Y}_{\mathrm{train}}.
\end{align}
Since the kernel $\boldsymbol{K}$ is positive semidefinite, let its spectral decomposition $\boldsymbol{K}=\boldsymbol{Q^{T}\Lambda Q}$, and since $e^{-\boldsymbol{K}t}=\boldsymbol{Q}^{\boldsymbol{T}}e^{-\boldsymbol{\Lambda}t}\boldsymbol{Q}$, we have
\begin{align}\label{eq:3.10}
\boldsymbol{Q}\left(f(\boldsymbol{X}_{\mathrm{train}},\boldsymbol{\theta}(t))-\boldsymbol{Y}_{\mathrm{train}}\right)=-e^{\boldsymbol{-\Lambda}t}\boldsymbol{Q}\boldsymbol{Y}_{\mathrm{train}}.
\end{align}
This means that components of an objective function corresponding to a kernel eigenvector with larger eigenvalues will learn more quickly. For traditional FNN, the eigenvalue of NTK decays rapidly. As a result, convergence to the high-frequency components of the objective function is very slow. Next, we use FNN to approximate the following function, and the results are shown in Figure \ref{fig1}. On the left, we can observe that even after 15,000 epochs, the network poorly fits the function with high-frequency components. From the figure on the right, we can see that the network first fits the low-frequency components of $f(x)$ and then gradually starts to fit the high-frequency components of $f(x)$.
\begin{align}\label{eq:3.11}
	f\left( x \right) =\sin \left( 2x \right) +\sin \left( 7x \right) +\sin \left( 16x \right).
\end{align}

\begin{figure}[htbp]
	\centering
	\includegraphics[scale=0.25]{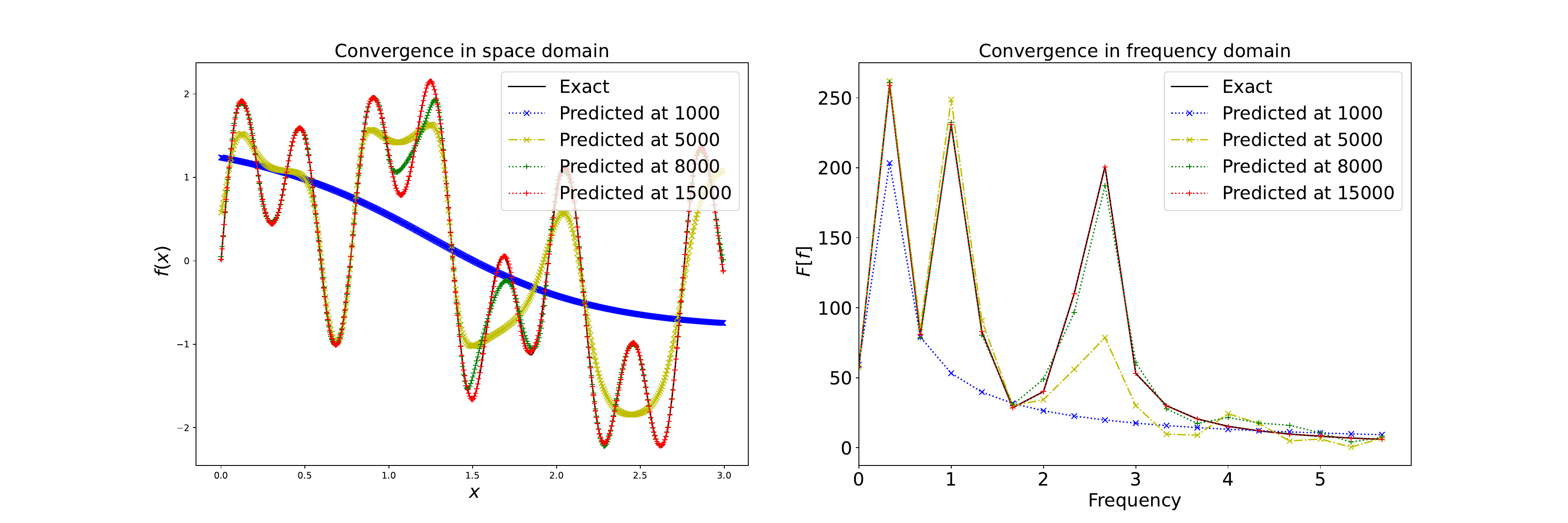}
	\caption{\; The convergence of FNN with Adaptive Activation Function in the spatial domain (left) and frequency domain (right) when approaching $f(x)$.}
	\label{fig2}
\end{figure}

The form of the Adaptive Activation Function\cite{jagtap2020adaptive, jagtap2022deep, jagtap2020locally} is as follows, and it has been utilized in PINNs to enhance the network's convergence speed. We employ an FNN with an Adaptive Activation Function to fit $f(x)$ under identical conditions, and the experimental results are depicted in Figure \ref{fig2}. Based on this figure, we can infer that the Adaptive Activation Function expedites the learning process of the high-frequency components of $f(x)$ and mitigates the network's spectral bias
\begin{align}\label{eq:3.12}
	\sigma(nax),
\end{align}
where $a$ is a trainable adaptive parameter, such an activation function is more likely to generate high-frequency nonlinear parts. This speeds up the network's learning of high-frequency information of the data. $n\geq1$ is the scaling factor selected by experience. For the hyper-parameter $n$, empirical studies showed $n=10$ optimally balances convergence speed and stability across our test cases, though adaptive tuning strategies may further improve performance.

\subsection{Proposed methods}\label{sec3.3}
As we previously discussed, WNO chooses to parameterize $R_{\phi}$ at the highest scale of the DWT. It retains only the information from the highest scale wavelet coefficients, allowing for a finite-dimensional parameterized space. This approach enables accurate data analysis and efficient learning. With the help of the wavelet transform's characteristics, WNO can effectively learn highly nonlinear differential operators, regardless of the geometry. However, WNO is not enough to extract high-frequency information from the data. 

\begin{figure}[htbp]
	\centering
	\includegraphics[scale=0.33]{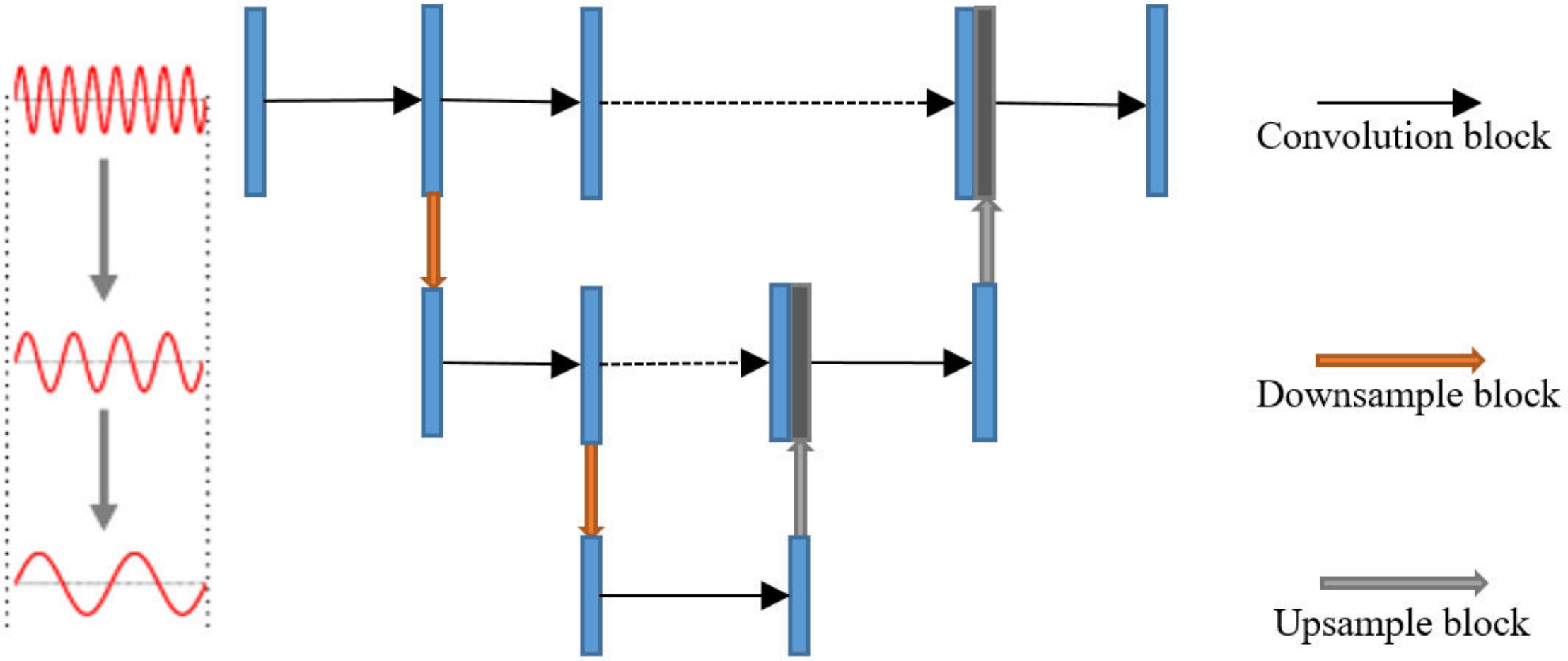}
	\caption{\ The architecture of U-Net path added to the U-WNO model.}
	\label{fig3}
\end{figure}

The U-Net is capable of handling local convolution, and it directly transmits feature maps from the encoder to the decoder, aiding in the preservation of more detailed information. Our introduction of the U-Net framework aims to improve WNO's ability to reconstruct high-frequency data. The U-Net structure is depicted in Figure \ref{fig3}, where it subsamples data to various scales. This allows for obtaining high-frequency information through convolution operations at different scales, thereby enhancing the expressiveness of neural operators on high-dimensional information.
At the same time, we also introduce residual shortcuts to ensure that each modified wavelet layer can enhance the performance of neural operators in the spatial domain. We use the Adaptive Activation Function described in Section \ref{sec3.2}. Thus, in the U-WNO layer, Eq.(\ref{eq:3.6}) is rewritten as
\begin{align}\label{eq:3.13}
	\begin{array}{c}
		v_{j+1}\left( x \right) :=\sigma \left(na \left( \mathcal{K}v_j \right) \left( x \right) +W\left( v_j\left( x \right) \right) \right) ,\ \ \ j=0\\
		v_{j+1}\left( x \right) :=\sigma \left(na \left( \mathcal{K}\left( v_j+v_0 \right) \right) \left( x \right) +W\left( v_j\left( x \right) \right) +\mathcal{U}\left( v_j\left( x \right) \right) \right) ,\ \ j=1,2,\cdots ,L\\
	\end{array}\ ,\forall x\in D.
\end{align}
Where $\mathcal{K}$, $W$ are the same as in WNO. $\sigma$ is an activation function. $\mathcal{U}$ is the U-Net convolutional neural network operator. The size of the convolution kernel used in the U-Net path is different from that in $W$, which helps to improve the ability of the network to extract features.

\begin{figure}[!htbp]
	\centering
	\includegraphics[scale=0.3]{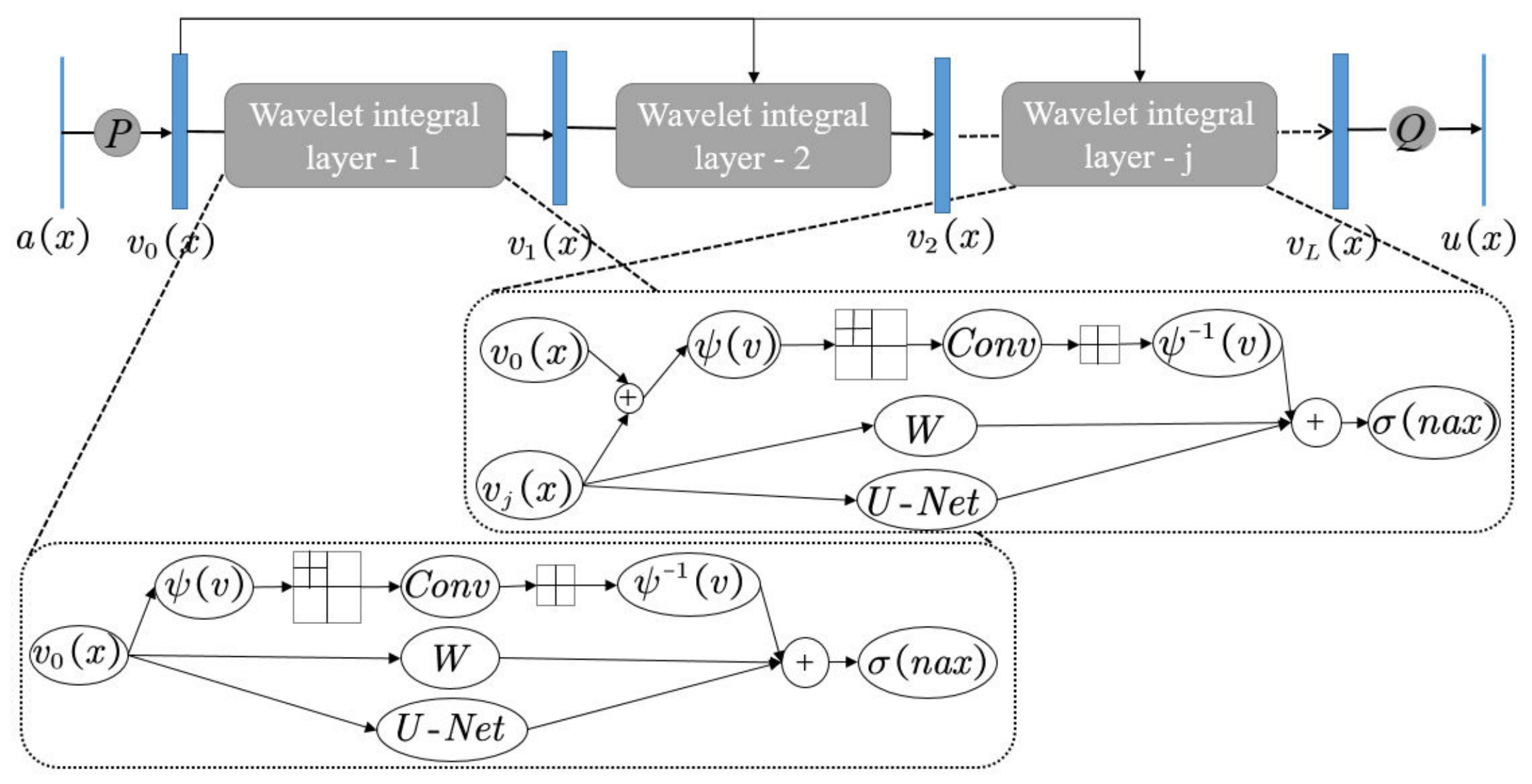}
	\caption{\ U-WNO model architecture. $a(x)$ is the input function, $P$ and $Q$ are FNNs, and $u(x)$ is the output function. In the U-WNO layer, we use the Wavelet transform denoted as $\psi$ and $R$ as a parameterization in Wavelet space. The inverse Wavelet transform is denoted as $\psi^{-1}$, and the linear bias term is referred to as $W$. The activation function is represented by $\sigma$. In the U-WNO layer, we utilize the residual shortcut to obtain complete information from $v_{0}$. Additionally, we add the U-Net path to enhance the performance of U-WNO in capturing high-frequency information.}
	\label{fig4}
\end{figure}

The U-WNO architecture also consists of three steps. First, as with WNO, the input function $a(x)$ is transformed to the higher-dimensional space $v_0 (x) =P (a(x))$, and then the U-WNO layer is applied: $v_{0}\mapsto v_{1}\mapsto...\mapsto v_{L}$, where $v_{j}$ for $j = 0,1,...,L$. Finally, $v_{L}$ is projected back to the original space using $Q$, resulting in $u(x)=Q(v_{L}(x))$. Figure \ref{fig4} provides a schematic of the U-WNO architecture. The step-by-step implementation of U-WNO is provided in Algorithm \ref{Algorithm 1}.

\begin{algorithm}[!htbp]
	\caption{Algorithm of the U-WNO}\label{Algorithm 1}
	\leftline{\textbf{Requirements:} $N$-samples of the pair $\left\{ a ( x ) \in \mathbb { R } ^ { n _ { D } \times d _ { a } } , u ( x ) \in \mathbb { R } ^ { n _ { D } \times d _ { u } } \right\}$, coordinates $x \in D $, and network}\par \hspace{2.6cm} hyperparameters.
	\begin{algorithmic}[1]
		\State Stack the inputs: $\{ a ( x ) , x \} \in \mathbb { R } ^ { n _ { D } \times 2 d_a } $.
		\For {epoch $= 1 , \ldots $, epochs}
			\State Uplift the input using transformation $P ( \cdot ) : v _ { 0 } ( x ) \in \mathbb { R } ^ { n_{D} \times d _ { v } } = P \left( \{ a ( x ) , x \} \in \mathbb { R } ^ { n _ { D } \times 2 d _ { a } } \right)$.
			\For {$j = 0 , \ldots , L$ perform the iterations: $v _ { j + 1 } = G \left( v _ { j } \right)$}
			\If {$j = 0$}
				\State Decompose the input using wavelet decomposition: $\mathcal { W } \left( v _ { j } ( x ) \right) \in \mathbb { R } ^ { n _ { D } / 2 ^ { m } \times d _ { v } } .$ 
				\State Parameterize the kernel $k _ { \phi }$ in the wavelet space : $R _ { \phi } * \mathcal { W } \left( v _ { j } ( x ) \right)$, see Eq.(\ref{eq:3.5}) 
				\State Reconstruct the convoluted input: $v _ { j + 1 } ^ { 1 } ( x ) = \mathcal { W } ^ { - 1 } \left( R _ { \phi } * \mathcal { W } \left( v _ { j } ( x ) \right) \right) ,$ see Eq.(\ref{eq:3.4})
				\State  Perform the linear transform:$v_{j+1}^2(x)=W*v_j(x)$ using Convolutional Neural Networks.
				\State Computes the output of the U-Net: $v_{j+1}^3(x)=\mathcal{U}\left( v_j\left( x \right) \right)$ . 
				\State  Add the outputs of step 8, 9 and 10: $\tilde { v } _ { j + 1 } ( x ) \in \mathbb { R } ^ { n _ { D } \times d _ { v } } = \left( v _ { j + 1 } ^ { 1 } + v _ { j + 1 } ^ { 2 } + v _ { j + 1 } ^ { 3 } \right) ( x ),$ see \par 
				\hspace{8.5mm}Eq.(\ref{eq:3.13})
			\Else
				\State Introduce residual shortcuts and proceed to Steps 6, 7, and 8: $\mathcal { W } \left( v _ { j } ( x ) + v _ { 0 } ( x ) \right) \in \mathbb { R } ^ { n _ { D } / 2 ^ { m } \times d _ { v } } .$ 
				\State Proceed to Step 7 and 8: $v _ { j + 1 } ^ { 1 } ( x ) = \mathcal { W } ^ { - 1 } \left( R _ { \phi } * \mathcal { W } \left( v _ { j } ( x ) + v _ { 0 } ( x ) \right) \right) .$
				\State  Proceed to Step 9 and 10: $v_{j+1}^2(x)=W*v_j(x)$, $v_{j+1}^3(x)=\mathcal{U}\left( v_j\left( x \right) \right)$. 
				\State  Add the outputs of step 14 and 15: $\tilde { v } _ { j + 1 } ( x ) \in \mathbb { R } ^ { n _ { D } \times d _ { v } } = \left( v _ { j + 1 } ^ { 1 } + v _ { j + 1 } ^ { 2 } + v _ { j + 1 } ^ { 3 } \right) ( x ) .$ 
				\If {$j \neq L$}
					 \State Apply the activation to complete the iteration: $v _ { j + 1 } \in \mathbb { R } ^ { n _ { D } \times d _ { v } } = \sigma \left(na \tilde { v } _ { j + 1 } ( x ) \right) .$ 
				\EndIf
			\EndIf
				\EndFor
				\State Apply an activated lifting transformation: $r(x)\in\mathbb{R}^{n_{D}\times d_{r}}=\sigma\left(Q_{1}(v_{L}(x))\right)$($Q_{1}(\cdot):\mathbb{R}^{d_{v}}\mapsto\mathbb{R}^{d_{r}}$is a \par  
				\hspace{-1.5mm}(FNN).
				\State Compute the final output: $\hat{u}(x)\in\mathbb{R}^{n_{D}\times d_{u}}=Q_{2}\left(r(x)\right)$($Q_{2}(\cdot):\mathbb{R}^{d_{r}}\mapsto\mathbb{R}^{d_{u}}$ is a FNN). 
				\State Compute the loss: $\mathcal { L } ( u , \hat { u } ) $. 
				\State Compute the gradient of the loss: $\frac { \partial \mathcal { L } ( u , \hat { u } ) } { \partial \theta _ { N N } } $. 
				\State Update the parameters of the network using the gradient.  
			\EndFor
	\end{algorithmic}
	{\textbf{Output:} Prediction of the solution of the PDEs}
\end{algorithm}

\section{Numerical experiments}
\vspace{-2mm}
In this section, the proposed U-WNO is applied to parametric PDEs, which include various problems in fluid dynamics, gas dynamics, traffic flow, atmospheric science, and phase field modeling. 

We compare the results obtained by four popular neural operators and their variants: (a) DeepONet, (b) FNO, (c) Multiwavelet transform operator (MWT) \cite{gupta2021multiwavelet}, (d) WNO. The U-WNO architecture consists of four U-WNO layers, each of which is activated by the GELU or MISH activation function. The mother wavelet is taken from the Daubechies family \cite{daubechies1992ten}. The dataset size, wavelet, wavelet tree level for all cases and the number of layers in the network against each equation are provided in Table \ref{tab1}. For training the U-WNO, we use the Adam optimizer with an initial learning rate of 0.001. During optimization, the learning rate decays at a rate of 0.5 every 50 epochs. The number of epochs used to train the U-WNO architecture is shown in Table \ref{tab2}. The batch size in the data loader is 20. Numerical experiments are performed on two NVIDIA Tesla V100-PCIE-32GB GPUs. In all experiments in this paper, the former part of the dataset is selected for the training set, and the latter part of the dataset is selected for the test set, which is consistent with the benchmark methods such as FNO and WNO. The results of solving PDEs with different parameters are summarized in Table \ref{tab3}. The results show that the mean $L_2$ relative test error of the proposed U-WNO ranges from 0.043$\%$ to 1.56$\%$, and U-WNO produces the lowest error in seven out of nine cases. In the remaining two problems, U-WNO is not the most accurate, but it is not far from the best result. Compared to WNO, U-WNO demonstrates significantly improved accuracy in all nine cases. Table \ref{tab4} shows the number of parameters and the time per epoch for WNO and U-WNO. Although the number of parameters of U-WNO is increased by about 75\% on average compared with WNO, the average increase of the training time of a single epoch of U-WNO is less than 40\%, and the final accuracy is significantly improved, which indicates that the increase in computational cost is acceptable.
The loss function used in U-WNO training and testing is the $L_{2}$ relative error function:
\begin{align}\label{eq:4.1}
	Loss=\frac{1}{N}\sum_{i=1}^N{\frac{\lVert u-\hat{u} \rVert _2}{\lVert u \rVert _2}}
\end{align}
where $N$ is the number of samples. $\hat{u}$ is the predicted value obtained by U-WNO inference, $u$ is the true value of the sample.
\begin{table}[ht!]
	\centering
	\caption{\; The relevant parameter information of the U-WNO.}
	\label{tab1}
	\begin{tabular}{llccccccc}
		\hline
		\multicolumn{2}{c}{\multirow{2}{*}{Examples}} & \multicolumn{2}{c}{Number of samples} & \multirow{2}{*}{Wavelet} & \multirow{2}{*}{m} & \multicolumn{2}{c}{Network dimensions} & \multirow{2}{*}{$\sigma \left( \cdot \right)$} \\ \cline{3-4} \cline{7-8}
		\multicolumn{2}{c}{}                          & Training           & Testing          &                          &                    & FNN1               & FNN2              &                                                \\ \hline
		\multirow{2}{*}{Burgers}   & time-independent & 1000               & 100              & db6                      & 8                  & 64                 & 128               & MISH                                           \\
		& time-dependent   & 480                & 20               & db6                      & 6                  & 40                 & 128               & MISH                                           \\ \hline
		\multirow{2}{*}{Advection} & time-independent & 900                & 100              & db6                      & 3                  & 96                 & 128               & MISH                                           \\
		& time-dependent   & 1000               & 100              & db6                      & 3                  & 80                 & 128               & MISH                                           \\ \hline
		\multirow{2}{*}{Darcy}     & rectangular      & 1000               & 100              & db4                      & 4                  & 64                 & 128               & GELU                                           \\
		& triangular       & 1000               & 100              & db6                      & 3                  & 64                 & 128               & MISH                                           \\ \hline
		\multicolumn{2}{c}{Allen-Cahn}                & 1000               & 100              & db6                      & 2                  & 64                 & 128               & MISH                                           \\ \hline
		\multicolumn{2}{c}{Non-homogeneous Poisson}   & 1000               & 100              & db4                      & 4                  & 64                 & 128               & GELU                                           \\ \hline
		\multicolumn{2}{c}{Navier-Stokes}             & 1000               & 100              & db4                      & 3                  & 26                 & 128               & GELU                                           \\ \hline
	\end{tabular}
\end{table}

\begin{table}[ht!]
	\centering
	\caption{\; The number of training epochs for different network structures and examples.}
	\label{tab2}
	\begin{tabular}{llccccc}
		\hline
		\multicolumn{2}{c}{\multirow{2}{*}{Examples}} & \multicolumn{5}{c}{Network architecture} \\ \cline{3-7} 
		\multicolumn{2}{l}{}                          & DeepONet   & FNO  & MWT  & WNO  & U-WNO  \\ \hline
		\multirow{2}{*}{Burgers}   & time-independent & 500000     & 500  & 500  & 500  & 500    \\
		& time-dependent   & -          & -  & -  & 500  & 500    \\ \hline
		\multirow{2}{*}{Advection} & time-independent & -          & -  & -  & 500  & 500    \\
		& time-dependent   & 250000     & 500  & 500  & 500  & 500    \\ \hline
		\multirow{2}{*}{Darcy}     & rectangular      & 100000     & 500  & 500  & 500  & 500    \\
		& triangular       & 20000      & 500  & 500  & 500  & 500    \\ \hline
		\multicolumn{2}{c}{Allen-Cahn}                & 100000     & 500  & 500  & 500  & 500    \\ \hline
		\multicolumn{2}{c}{Non-homogeneous Poisson}   & -          & 500  & -  & 500  & 500    \\ \hline
		\multicolumn{2}{c}{Navier-Stokes}             & 100000     & 500  & 500  & 500  & 500    \\ \hline
	\end{tabular}
\end{table}

\begin{table}[ht!]
	\centering
	\caption{\; Mean $L_{2}$ relative error between the exact and predicted results.}
	\label{tab3}
	\begin{threeparttable}
		\begin{tabular}{llcccccc}
			\hline
			\multicolumn{2}{c}{\multirow{2}{*}{Examples}} & \multicolumn{6}{c}{Network architecture}                        \\ \cline{3-8} 
			\multicolumn{2}{c}{}                          & DeepONet & FNO & \multicolumn{1}{l}{dgFNO+\tnote{a}} & MWT & WNO & U-WNO \\ \hline
			\multirow{2}{*}{Burgers}   & time-independent & 2.15$\%$ & 1.60$\%$ & -  & 0.19$\%$ & 1.87$\%$ & 1.56$\%$   \\
			& time-dependent   & -        & -        & -    & - & 0.23$\%$ & 0.043$\%$   \\ \hline
			\multirow{2}{*}{Advection} & time-independent & -        & - & - & - & 0.85$\%$ & 0.053$\%$   \\
			& time-dependent   & 0.32$\%$ & 47.7$\%$ & 0.62$\%$ & 10.22$\%$ & 0.62$\%$ & 0.059$\%$   \\ \hline
			\multirow{2}{*}{Darcy}     & rectangular      & 2.98$\%$   & 1.08$\%$ & -  & 0.89$\%$ & 1.84$\%$ & 0.62$\%$   \\
			& triangular       & 2.64$\%$ & - & 7.82$\%$ & 0.87$\%$ & 0.77$\%$ & 0.20$\%$   \\ \hline
			\multicolumn{2}{c}{Allen-Cahn}                & 17.7$\%$   & 0.93$\%$ & - & 4.84$\%$ & 0.22$\%$ & 0.095$\%$   \\ \hline
			\multicolumn{2}{c}{Non-homogeneous Poisson}   & -        & 0.16$\%$ & -                        & - & 0.36$\%$ & 0.044$\%$   \\ \hline
			\multicolumn{2}{c}{Navier-Stokes}             & 1.78$\%$   & 1.28$\%$ & - & 0.63$\%$ & 3.16$\%$ & 0.96$\%$   \\ \hline
		\end{tabular}
		\begin{tablenotes}
			\item [a] Note that in complex geometric conditions, FNO does not work, thus we use dgFNO+\cite{lu2022comprehensive}.
		\end{tablenotes}
	\end{threeparttable}
\end{table}

\begin{table}[ht!]  
\centering  
	\caption{\; The number of parameters for different network structures and the time required for each epoch.}  
	\label{tab4}  
    \begin{threeparttable}  
\begin{tabular}{llcccc}  
\hline  
\multicolumn{2}{c}{\multirow{2}{*}{Examples}} & \multicolumn{2}{c}{WNO}  & \multicolumn{2}{c}{U-WNO}   \\ \cline{3-6}   
\multicolumn{2}{c}{}                          & Parameters & Time & Parameters & Time \\ \hline  
\multirow{2}{*}{Burgers}   &time-independent & 484033    & 6.32s     & 1159362    & 7.30s         \\
                           &time-dependent  & 243217   & 32.97s      & 506258    & 44.29s       \\ \hline  
\multirow{2}{*}{Advection} & time-independent & 1082273  & 3.73s         & 2599074   & 5.03s         \\
                           & time-dependent & 756497    & 63.68s         & 1807378  & 96.71s         \\ \hline  
\multirow{2}{*}{Darcy}     & rectangular      & 7955201   & 6.74s        &10596866   & 9.54s         \\
                           & triangular       & 16802561  & 11.20s        & 18738434  & 16.63s         \\ \hline  
\multicolumn{2}{c}{Allen-Cahn}               & 16802561  & 6.24s         & 19261442 & 10.04s          \\ \hline  
\multicolumn{2}{c}{Non-homogeneous Poisson}   & 14770945   & 25.51s       & 17412610   & 39.33s        \\ \hline  
\multicolumn{2}{c}{Navier-Stokes}             & 1448207   & 33.72s        & 2562532   &42.70s        \\ \hline  
\end{tabular}  
\end{threeparttable}  
\end{table}

\subsection{1D Burgers equation}\label{sec4.1}
For the first example, we choose the Burgers equation with different initial conditions. The 1D Burgers equation is frequently employed to mathematically represent fluid dynamics and gas dynamics, traffic flow, and other braided formations. The 1D Burgers equation with periodic boundary conditions can be expressed mathematically as follows:

\begin{align}\label{eq:4.2}
	\begin{aligned}
		\partial_{t}u(x,t)+\dfrac{1}{2}\partial_{x}u^{2}(x,t)& =v\partial_{xx}u(x,t),\quad x\in(0,1),t\in(0,1] \\
		u(x=0,t)& =u(x=1,t),\quad x\in(0,1),t\in(0,1] \\
		u(x,0)& =u_0(x),\quad\quad\quad  x\in(0,1). 
	\end{aligned}
\end{align}
Where $\nu\in\mathbb{R}_+$ is the viscosity of the flow and $u_{0}\in L_{\mathrm{ner}}^{2}((0,1);\mathbb{R})$ is the initial condition. The initial conditions $u_{0}(x)$ are generated using a Gaussian random field: $\mu=\mathcal{N}(0,625(-\Delta+25I)^{-2})$. Consider periodic boundary conditions of the form $u(x-\pi,t)=u(x+\pi,t);x \in (0,1),t \in (0,1]$. The aim here is to learn the operator $\mathcal{D}:u_{0}(x)\mapsto u(x,1)$. We choose $\nu=0.1$ and a spatial resolution of 1024 for the experiments. The dataset is taken from Ref. \cite{li2020fourier}.

\begin{figure}[ht!]
	\centering
	\includegraphics[scale=0.5]{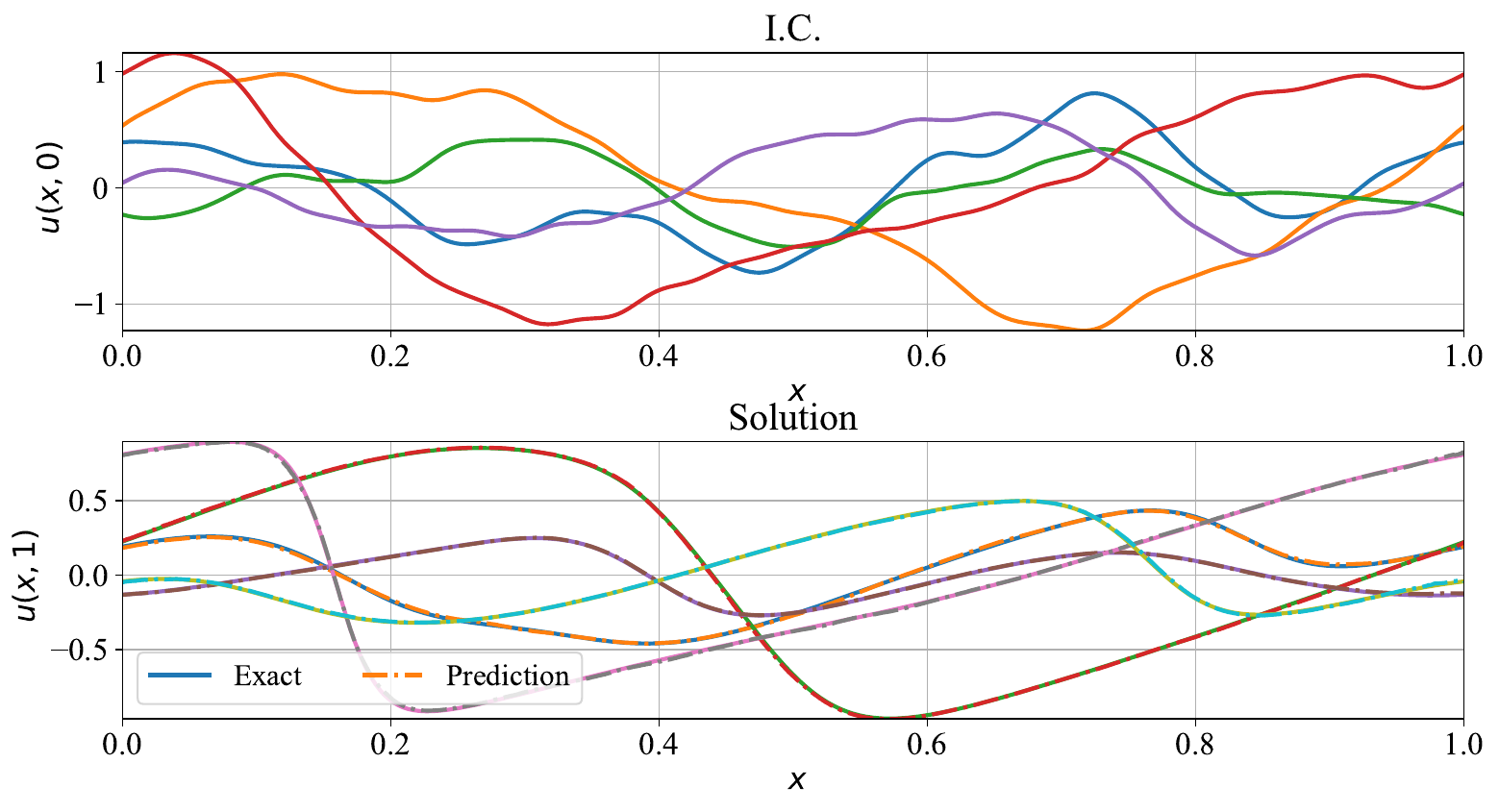}
	\caption{\;1D Burgers equation with periodic boundary conditions at 1024 spatial resolution. (top) Initial condition at $t = 0$, where different colored lines represent different initial value conditions. (bottom) Color-matched solid/dashed lines represent true/predicted solutions for each initial condition of 1D Burgers equation at $t = 1$.}
	\label{fig5}
\end{figure}

\textbf{Results}: The prediction results and the mean test $L_{2}$ relative error are provided in Figure \ref{fig5} and Table \ref{tab3}. From Table \ref{tab3}, the MWT obtains the lowest error, followed by the proposed U-WNO. The error of U-WNO is lower than DeepONet and FNO, and U-WNO significantly reduces the error compared with WNO. Although the error predicted using U-WNO is higher than that of the MWT from a quantitative point of view, Figure \ref{fig5} shows that the difference between the predicted solutions and true solutions is small. This shows that although U-WNO is not the best, it still has good expressiveness in this problem.

\textbf{1D time-dependent Burgers equation with discontinuity in the solution field:} For small values of the viscosity parameter $\nu$, the Burgers equation causes a sharp jump in the solution field. Modeling such sharp jumps in the solution field is extremely challenging. In this example, our objective is to understand how the operator's initial conditions lead to a sharp jump in the solution field. The 1D Burgers equation, in the following form, can be employed to simulate sharp jumps in the solution field:
\begin{align}\label{eq:4.3}
	\begin{aligned}
		\partial_{t}u(x,t)+u(x,t)\partial_{x}u(x,t)& =\frac{0.01}{\pi}\partial_{xx}u(x,t),\quad x\in[-1,1],t\in[0,1] \\
		u(x=-1,t)& =u(x=1,t)=0,\quad x\in[-1,1],t\in[0,1] \\
		u(x,0)& =-sin(\pi x)+\zeta sin(\pi x),\quad x\in[-1,1]. 
	\end{aligned}
\end{align}
Where $\zeta\in[0,0.5]$ is a random number. The aim is to learn the operator mapping the solution fields at the time steps $t\in[0,10]$ to the solution fields up to time steps $t\in(10,50]$ defined by $\mathcal{D}: u|_{(-1,1)\times[0,10]}\mapsto u|_{(-1,1)\times[10,50]}$. Data with a time step of $\Delta t=0.02$ and a spatial resolution of 512 were used for training and testing. The dataset taken from Ref. \cite{tripura2023wavelet}. 

\begin{figure}[ht!]
	\centering
	\includegraphics[scale=0.5]{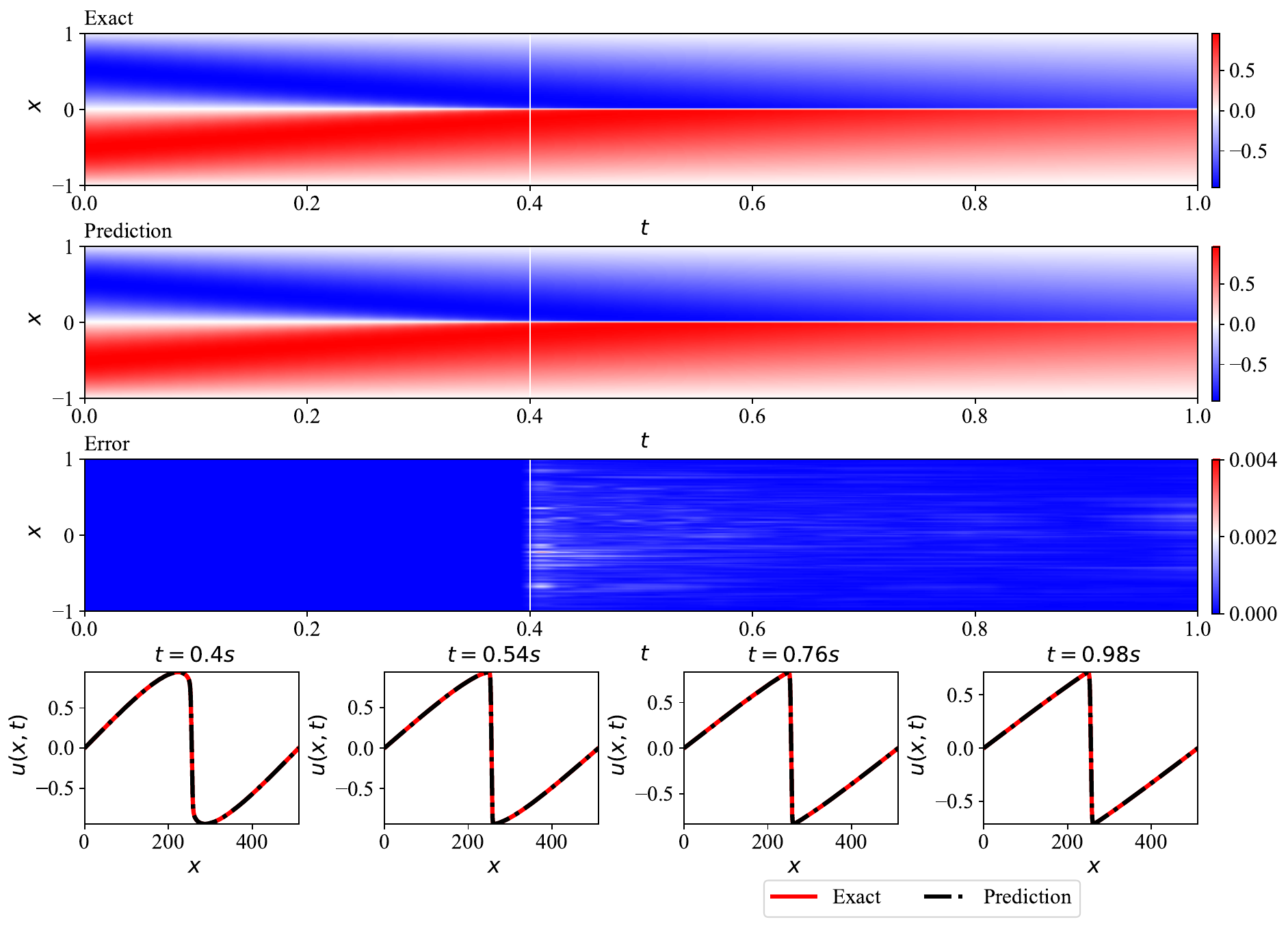}
	\caption{\;1D Burgers equation with a discontinuity in the solution field. The first three subplots illustrate the truth solutions, prediction ones, and errors, respectively. The four subplots in the last row depict the predicted and true solutions for the four moments with relatively large errors.}
	\label{fig6}
\end{figure}

\textbf{Results}: The true, predicted, and error results for an initial condition are shown in the first three subplots of Figure \ref{fig6}. From Figure \ref{fig6} and Table \ref{tab3}, we can conclude that the prediction results obtained using U-WNO are accurate. The third subplot shows that the error is relatively large after $t=0.4\ s$ (the white line in the figure), while the four plots in the last row depict the true and predicted solutions at $t=0.4\ s$, $t=0.54\ s$, $t=0.76\ s$ and $t=0.98\ s$. Visually, there is no difference between the predicted and true solutions. These show that the proposed U-WNO can accurately solve the 1D Burgers equation with sharp jump cases within the solution field.

\subsection{1D wave advection equation}\label{sec4.2}
The wave advection equation describes the solution to a scalar in a given velocity field. The wave advection equation with periodic boundary conditions has the following form:
\begin{align}\label{eq:4.4}
	\begin{aligned}
		\partial_{t}u(x,t)+\nu\partial_{x}u(x,t)=0,\quad x\in(0,1),t\in(0,1)\\
		u(x-\pi,t)=u(x+\pi,t),\quad x\in(0,1),t\in(0,1).
	\end{aligned}
\end{align}
Here $\nu\in\mathbb{R}^{+*}$ denotes the velocity of the stream. The initial condition is,
\begin{align}\label{eq:4.5}
u(x,0)=h1_{\left\{c-\frac{\omega}{2},c+\frac{\omega}{2}\right\}}+\sqrt{\max(h^{2}-(a(x-c))^{2},0)}.
\end{align}
Here, the value of $\{c,\omega,h\}$ is randomly chosen as $[0.3,0.7]\times[0.3,0.6]\times[1,2]$. For $\nu = 1$, the solution of the advection equation is $u(x,t)=u_{0}(x-t)$, the time step is 0.025, and the spatial discrete points are 40, resulting in a resolution of $40 \times 40$ advection equation solution. The dataset is taken from Ref. \cite{lu2022comprehensive}. We first learn the mapping from the initial condition $u(x,0)$ to the solution at $t = 1$:
$$\text{Case I:}\quad\mathcal{G}_1:u(x,0)\mapsto u(x,t=1),$$
Then learn the operator mapping from the initial condition $u(x,0)$ to the solution of the whole domain
$$\text{Case II:}\quad\mathcal{G}_2:u(x,0)\mapsto u(x,t).$$

\begin{figure}[ht!]
	\centering
	\includegraphics[scale=0.5]{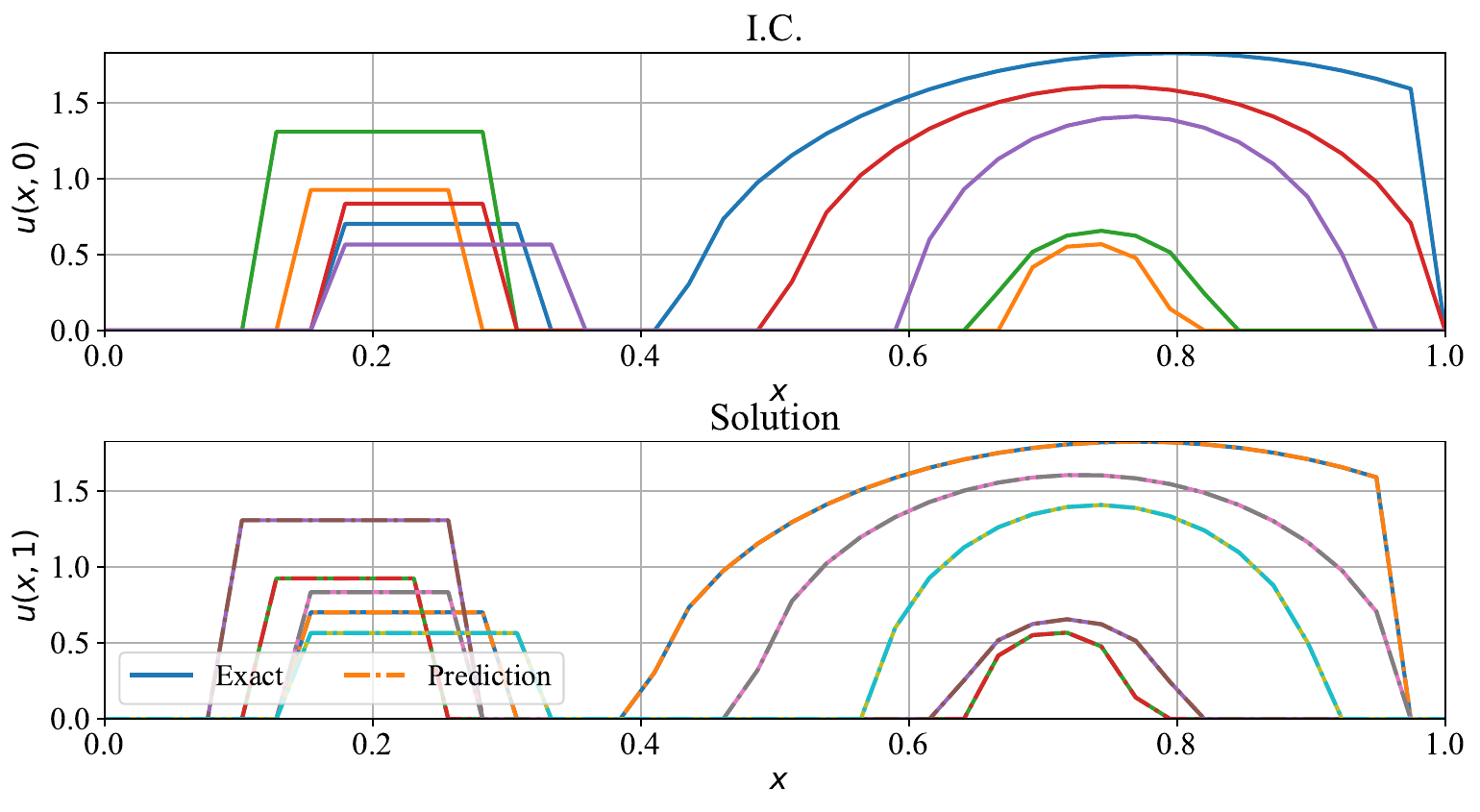}
	\caption{\;1D wave advection equation with periodic boundary conditions at $40\times 40$ spatial resolution. (top) Initial condition at $t = 0$. (bottom) True and predicted solutions of 1D wave advection equation at $t = 1$.}
	\label{fig7}
\end{figure}

\textbf{Results}: For case I, The prediction results and mean test $L_2$ relative error are shown in Figure \ref{fig7} and Table \ref{tab3}. From Table \ref{tab3}, we can conclude that U-WNO significantly reduces the error compared to WNO. Figure \ref{fig7} shows that there is almost no difference between the predicted and true results. This shows that U-WNO solves this problem well.

\begin{figure}[ht!]
	\centering
	\includegraphics[scale=0.35]{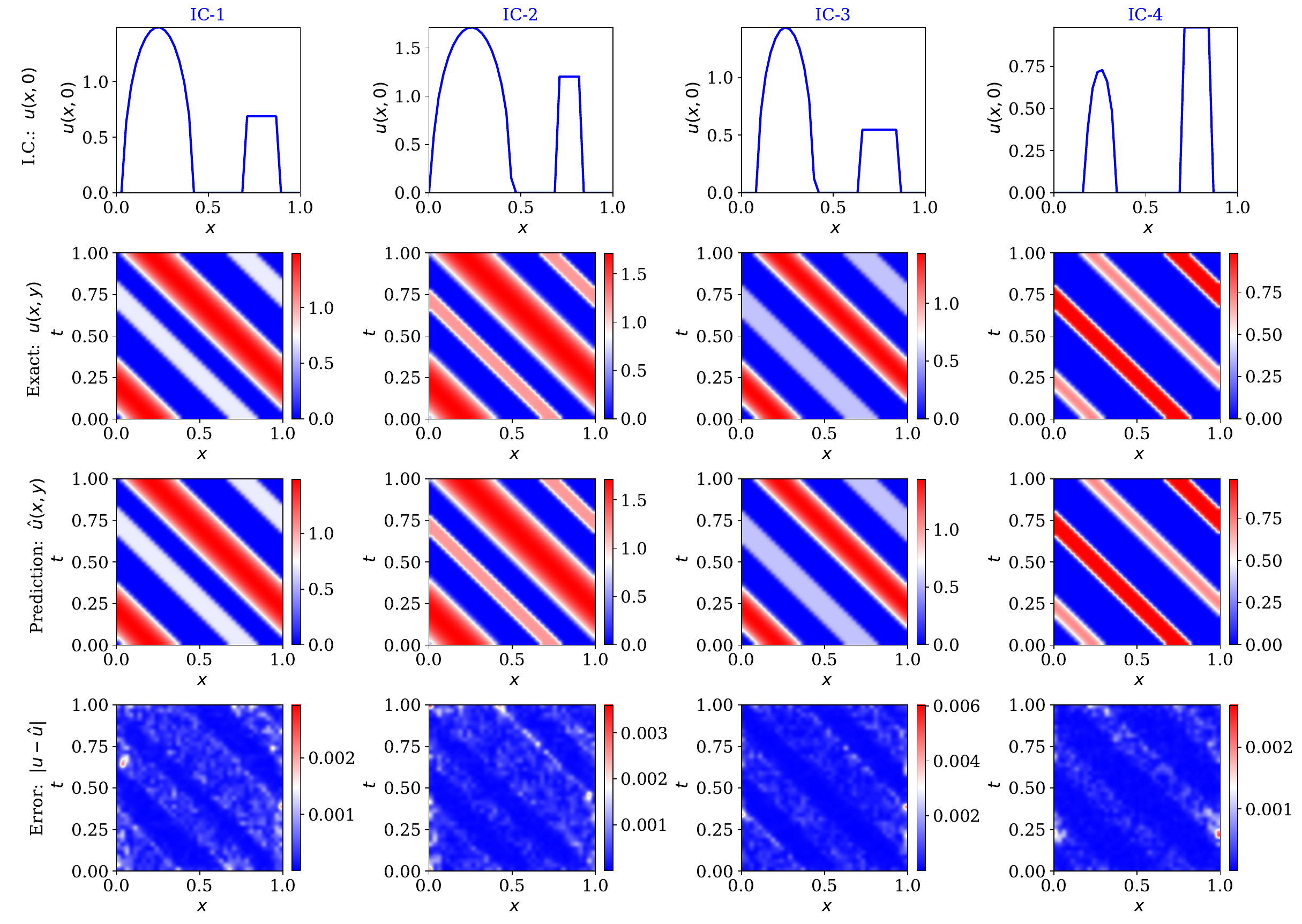}
	\caption{\; Wave advection equation with periodic boundary conditions with a resolution of $40 \times 40$. The figure shows the true solution, the predicted solution, and the error corresponding to the initial conditions in Eq.(\ref{eq:4.5})}.
	\label{fig8}
\end{figure}

For case II: The prediction results and mean $L_2$ relative errors are shown in Figure \ref{fig8} and Table \ref{tab3}. Except for FNO and MWT which have slightly higher errors, DeepONet,dgFNO+, WNO, and U-WNO all have errors less than 1$\%$, while the proposed U-WNO achieves the lowest error with an error less than 0.1$\%$.

\subsection{2D darcy flow equation}\label{sec4.3}
The 2D Darcy flow equation is widely used for fluid modeling. In this subsection, we consider Darcy flows of the following form
\begin{align}\label{eq:4.6}
	\begin{aligned}
		-\nabla\cdot(a(x,y)\nabla u(x,y))=f(x,y),\quad x,y\in(0,\mathbb{R})\\u(x,y)=u_{0}(x,y),\quad x,y\in\partial(0,\mathbb{R}).
	\end{aligned}
\end{align}
Where $a(x,y)$ represents the coefficients, the Darcy flow is defined in a rectangular field: $x\times y\in(0,1)^{2}$. In this example, we select the initial condition $u_{0}(x,y)=0$ and the source function $f(x,y)=1$. The aim is to learn the operator $D:a(x,y)\mapsto u(x,y)$.
The coefficients $a(x,y)$ are generated according to $a \sim \psi_{\#}\mathcal{N}(0,(-\Delta+9I)^{-2})$ with zero Neumann boundary conditions on the Laplacian. The map $\psi:\mathbb{R}\to\mathbb{R}$ takes the value $12$ in the positive part of the real line and $3$ in the negative part, defined pointwise in the forward direction. The dataset is taken from Ref. \cite{li2020fourier}. During training, the spatial resolution is $85\times 85$. 

\begin{figure}[ht!]
	\centering
	\includegraphics[scale=0.35]{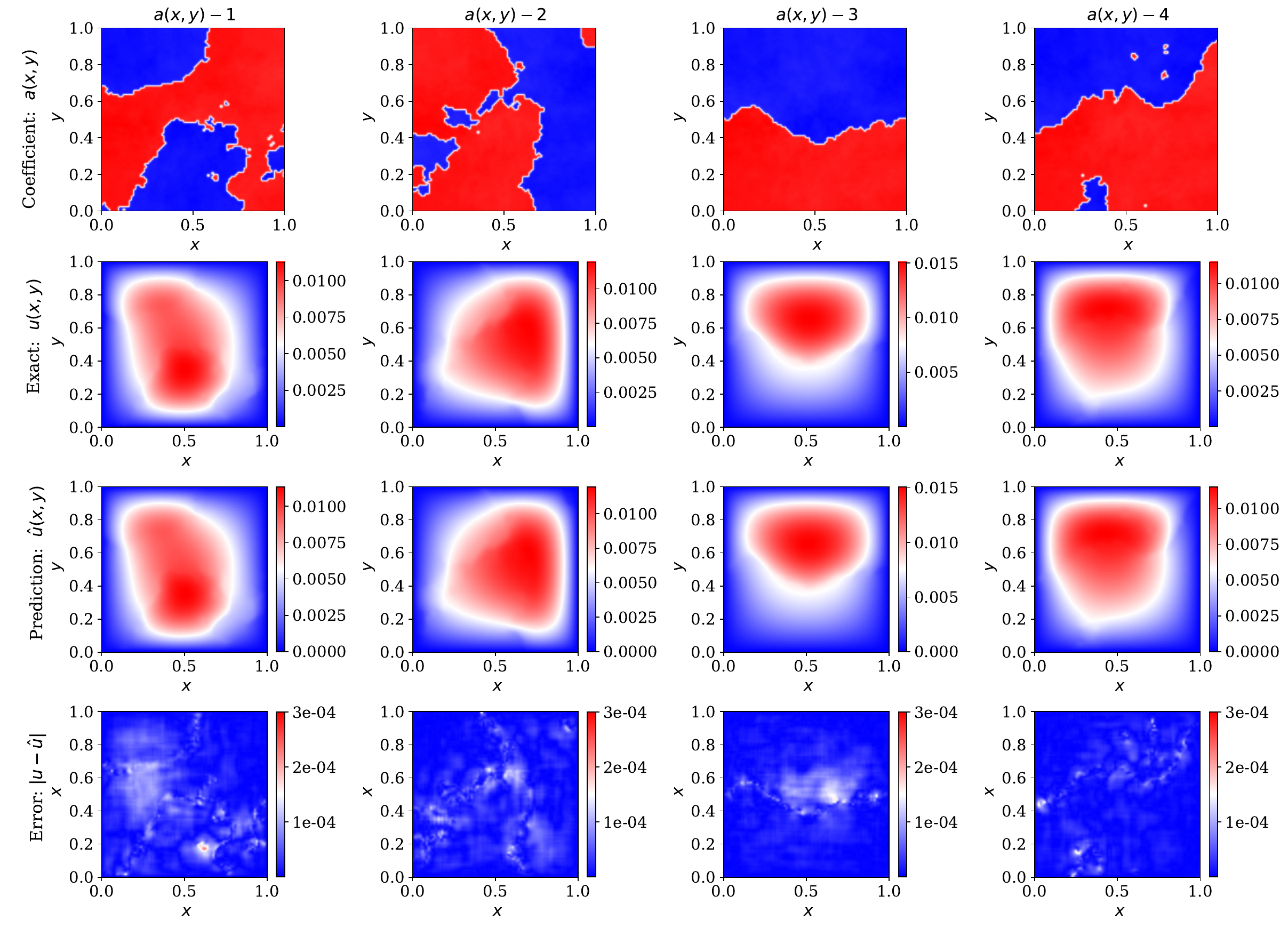}
	\caption{\; Darcy flow in a rectangular domain with spatial resolution $85\times 85$. For different coefficients, the solution field is obtained using the proposed U-WNO. The figure shows the predicted solutions, true solutions, and errors.}
	\label{fig9}
\end{figure}

\textbf{Results}: The prediction results of the 2D Darcy flow in the rectangular domain are shown in Figure \ref{fig9}, and the mean test $L_2$ relative error are shown in the Table  \ref{tab3}. In this case, it can be observed that the error of U-WNO is the lowest. In Figure \ref{fig9}, it can be concluded that the difference between the true and predicted solutions is small. 

\textbf{2D Darcy flow equation with a notch in triangular domain: }Here, we consider the Darcy flow equation with a notch in triangular domain with coefficients $a(x, y)=0.1$ and source function $f (x, y)=-1$. The aim is to learn the operator that maps the boundary conditions to the solution field, defined as $\mathcal{D}:u(x,y)|_{\partial\omega}\mapsto u(x,y)$. The dataset is taken from Ref. \cite{lu2022comprehensive}.

\begin{figure}[ht!]
	\centering
	\includegraphics[scale=0.35]{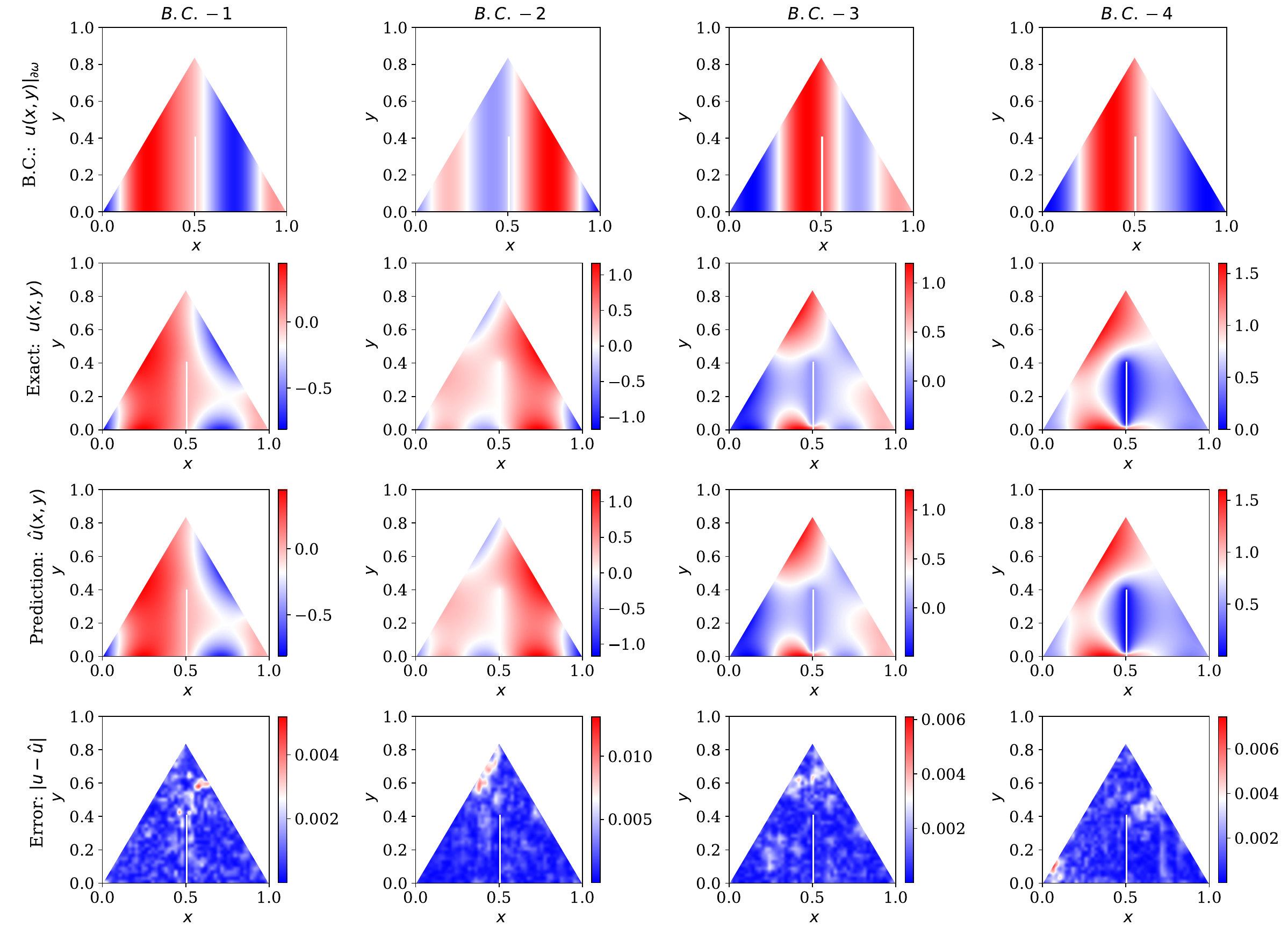}
	\caption{\; Darcy flow in a triangular domain with a spatial resolution of $51 \times 51$. By using the proposed U-WNO model, the predicted solutions, true solutions, and errors are obtained under different boundary conditions.}
	\label{fig10}
\end{figure}

\textbf{Results}: The prediction solutions and mean test $L_{2}$ relative error are shown in Figure \ref{fig10} and Table \ref{tab3}. From the values in Table \ref{tab3}, the error of U-WNO is the lowest, and the error is less than $1/3$ of the WNO. The error between the true value and the prediction can also be observed in Figure \ref{fig10}, and it can be further noted that the predicted solution yields accurate predictions both in the smooth region and near the notch, which indicates that the proposed U-WNO can learn nonlinear operators under complex geometric conditions with higher accuracy than WNO.

\subsection{2D Allen-Cahn equation}\label{sec4.4}
The Allen-Cahn equation is commonly used for predicting chemical reactions and phase separation in multicomponent alloys. The 2D Allen-Cahn equation is as follows.
\begin{align}\label{eq:4.7}
	\begin{matrix}
		\partial _tu\left( x,y,t \right) =\epsilon \varDelta u\left( x,y,t \right) +u\left( x,y,t \right) -u\left( x,y,t \right) ^3,&		x,y\in \left( 0,3 \right) ,t\in \left[ 0,10 \right]\\
		u\left( x=0,y,t \right) =u\left( x=3,y,t \right) ,&		y\in \left( 0,3 \right) ,t\in \left[ 0,10 \right]\\
		u\left( x,y=0,t \right) =u\left( x,y=3,t \right) ,&		x\in \left( 0,3 \right) ,t\in \left[ 0,10 \right]\\
		u\left( x,y,0 \right) =u_0\left( x,y \right) ,&		x,y\in \left( 0,3 \right)\\
	\end{matrix}
\end{align}
Where $\epsilon\in\mathbb{R}^{+}$ is responsible for the amount of diffusion. The problem is defined on a periodic boundary of $\epsilon=10^{-3}$, and the initial conditions are generated from a Gaussian random field with the following kernel:
\begin{align}\label{eq:4.8}
	\mathcal{K}(x,y)=\tau^{(\alpha-1)}\left(\pi^2(x^2+y^2)+\tau^2\right)^{\frac{\alpha}{2}},
\end{align}
where $\tau = 15$ and $\alpha = 1$. The aim here is to learn the operator $\mathcal{D} : u_{0}(x, y) \mapsto u(x, y, 10)$. Training and testing were performed on a resolution of $43 \times 43$. The dataset is taken from Ref. \cite{tripura2023wavelet}.

\begin{figure}[ht!]
	\centering
	\includegraphics[scale=0.35]{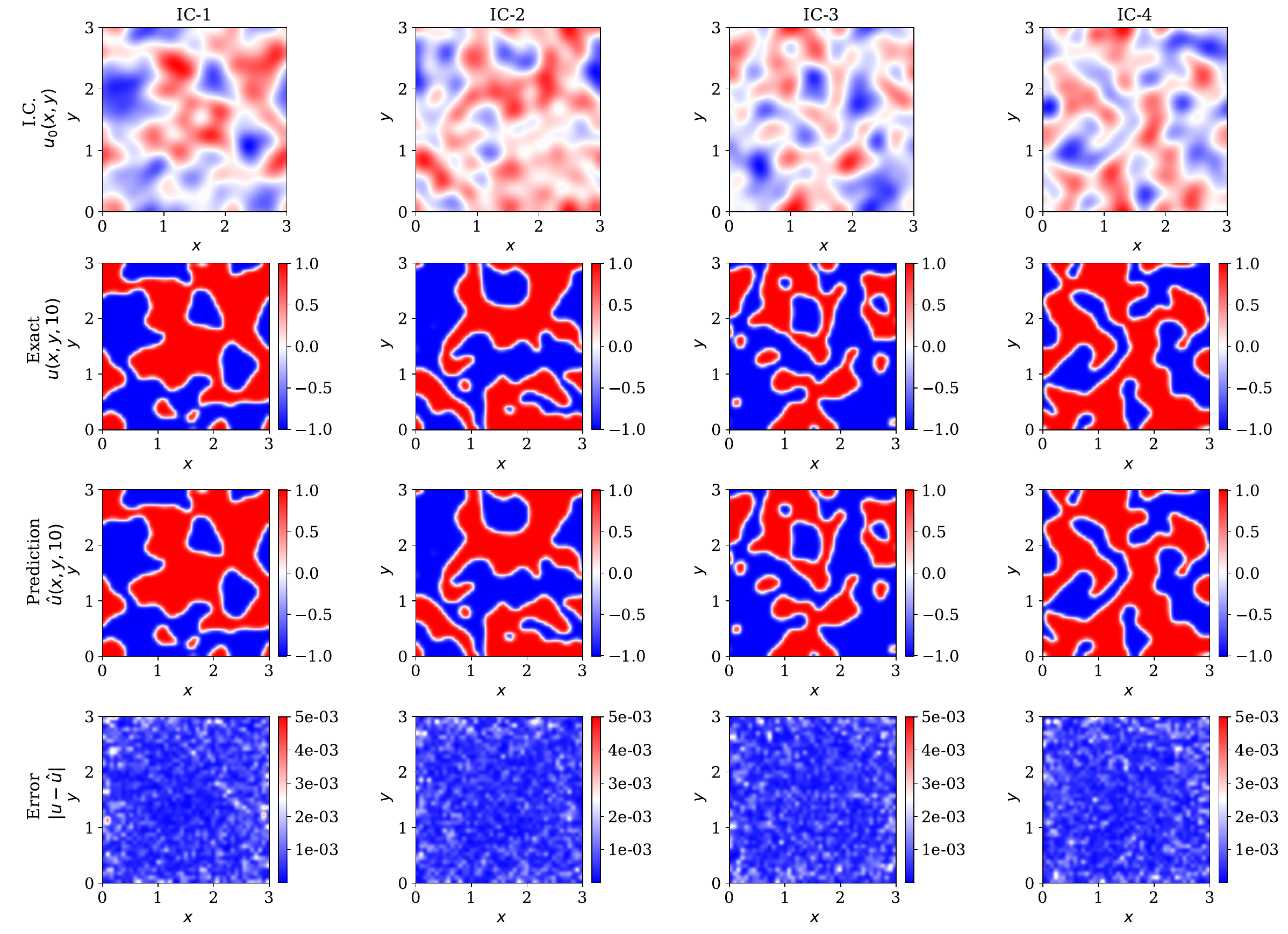}
	\caption{\; Allen-Cahn equation with periodic boundary conditions at resolution $43\times 43$. The predicted solutions, true solutions, and errors are obtained under different boundary conditions by using the proposed U-WNO model.}
	\label{fig11}
\end{figure}

\textbf{Results}: The mean test $L_{2}$ relative error and prediction results of the proposed U-WNO for different initial conditions are shown in Table \ref{tab3} and Figure \ref{fig11}. The values in Table\ref{tab3} show that U-WNO obtains the lowest prediction error, which is less than $1/2$ of the error obtained by WNO and slightly less than 0.1$\%$. In Figure \ref{fig11}, it can be concluded that the difference between the true and predicted results is small for all given initial conditions. This shows that the proposed U-WNO can accurately solve the nonlinear problem in 2D domain.

\subsection{Non-homogeneous Poisson equation}\label{sec4.5}
The Poisson equation with a source term $f (x, y)$ and periodic boundary condition is expressed as:
\begin{align}\label{eq:4.9}
	\begin{aligned}
		&\partial_{xx}u+\partial_{yy}u=f(x,y),\quad x\in[-1,1],\quad y\in[-1,1]\\&u(x=-1,y)=u(x=1,y)=u(x,y=-1)=u(x,y=1)=0.
	\end{aligned}
\end{align}
In this example, the aim is to learn that the mapping between the source function $f (x, y)$ and the solution $u(x, y)$ is defined as $D: f(x,y)\mapsto u(x,y)$, Source function we choose $f(x,y) = 16\beta\pi^{2} (\cos(4\pi y)\sin(4\pi x)+\sin(4\pi x)(\cos(4\pi y)-1))-\alpha\pi^{2}(\cos(\pi y)\sin(\pi x)+\sin(\pi x)(\cos(\pi y)+1))$. Then the solution is $u(x,y)=\alpha\sin(\pi x)(1+\cos(\pi y))+\beta\sin(2\pi x)(1-\cos(2\pi y))$. $\alpha ,\ \beta \in \left[ -2,2 \right] $ are uniformly distributed random numbers. The samples of the training and testing instances are yielded by varying the variable $\alpha$ and $\beta$. A resolution of $85\times 85$ was used for training and testing. 

\textbf{Results}: The results in Figure \ref{fig11} include the source function, the true solution, the predicted solution, and the prediction error. The results clearly show that the predicted and true solutions obtained by the proposed U-WNO are in good agreement. In addition, from the results shown in Table \ref{tab3}, the mean $L_2$ relative test error of U-WNO is lower than the error of WNO and FNO, reaching 0.044$\%$, which indicates that U-WNO can learn nonlinear operators from source functions to solutions in 2D domain.

\begin{figure}[ht!]
	\centering
	\includegraphics[scale=0.35]{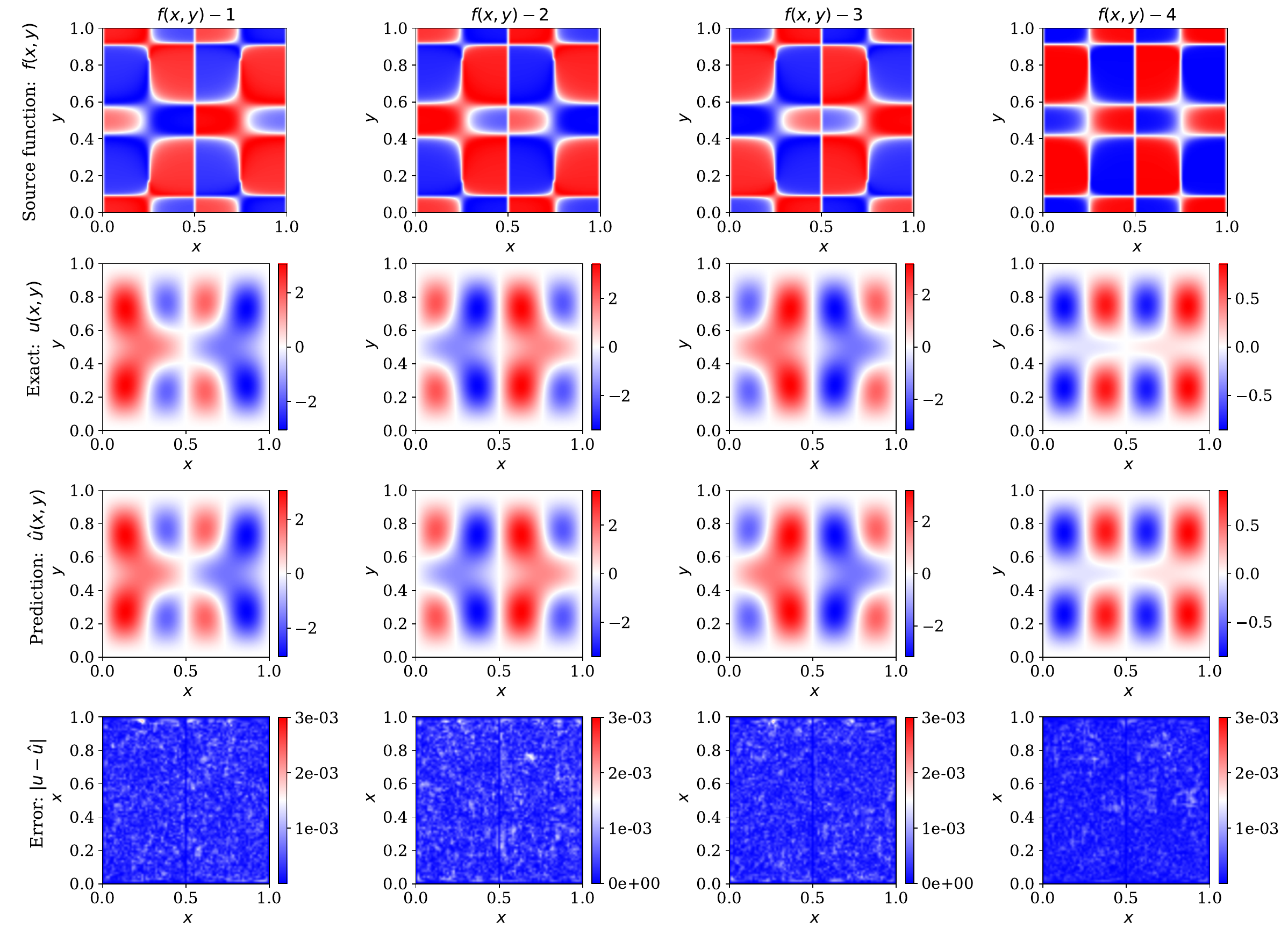}
	\caption{\; The proposed U-WNO maps the source function $f (x, y)$ of the Non-homogeneous Poisson equation to the corresponding solution $u(x, y)$ with a spatial resolution of $85 \times 85$, the corresponding real solution, prediction, and error composition of the result.}
	\label{fig12}
\end{figure}

\subsection{2D Navier–Stokes equation}\label{sec4.6}
Navier-Stokes equation is often used in air flows around aircraft wings, ocean currents, thermodynamics, etc. The 2D incompressible Navier-Stokes equations in the vorticity velocity form are of the following form.
\begin{align}\label{eq:4.10}	\begin{aligned}\partial_{t}\omega(x,y,t)+u(x,y,t)\cdot\nabla\omega(x,y,t)&=\nu\Delta\omega(x,y,t)+f(x,y),&x,y\in(0,1),t\in(0,T]\\\nabla\cdot u(x,y,t)&=0,&x,y\in(0,1),t\in[0,T]\\\omega(x,y,0)&=\omega_{0}(x,y),&x,y\in(0,1)\end{aligned}
\end{align}
where $\nu\in\mathbb{R}$, $f (x, y)$ is the fluid viscosity and the source function. $u(x, y, t)$ and $\omega(x,y,t)$ are the velocity and vorticity fields of the fluid, and the viscosity is taken as $\nu = 10^{-3}$. The aim is to learn the operator mapping time step $t \in [0, 10]$ of vorticity field  to the time step $t \in (10, 20]$ of the vorticity field, defined as $D: \omega|_{(0,1)^{2}\times[0,10]}\mapsto\omega|_{(0,1)^{2}\times[10,20]}$. The initial vorticity $\omega_{0}(x,y)$ by Gaussian random field simulation for the $\omega_{0}(x,y)=\mathcal{N}\left(0,7^{3/2}(-\Delta+49I)^{-2.5}\right)$. The source function $f(x,y)=0.1 (\sin{(2\pi (x+y))}+\cos{(2\pi (x+y))})$. The dataset is taken from Ref. \cite{tripura2023wavelet}. For training and testing, the spatial resolution is fixed to $64 \times 64$.

\begin{figure}[ht!]
	\centering
	\includegraphics[scale=0.35]{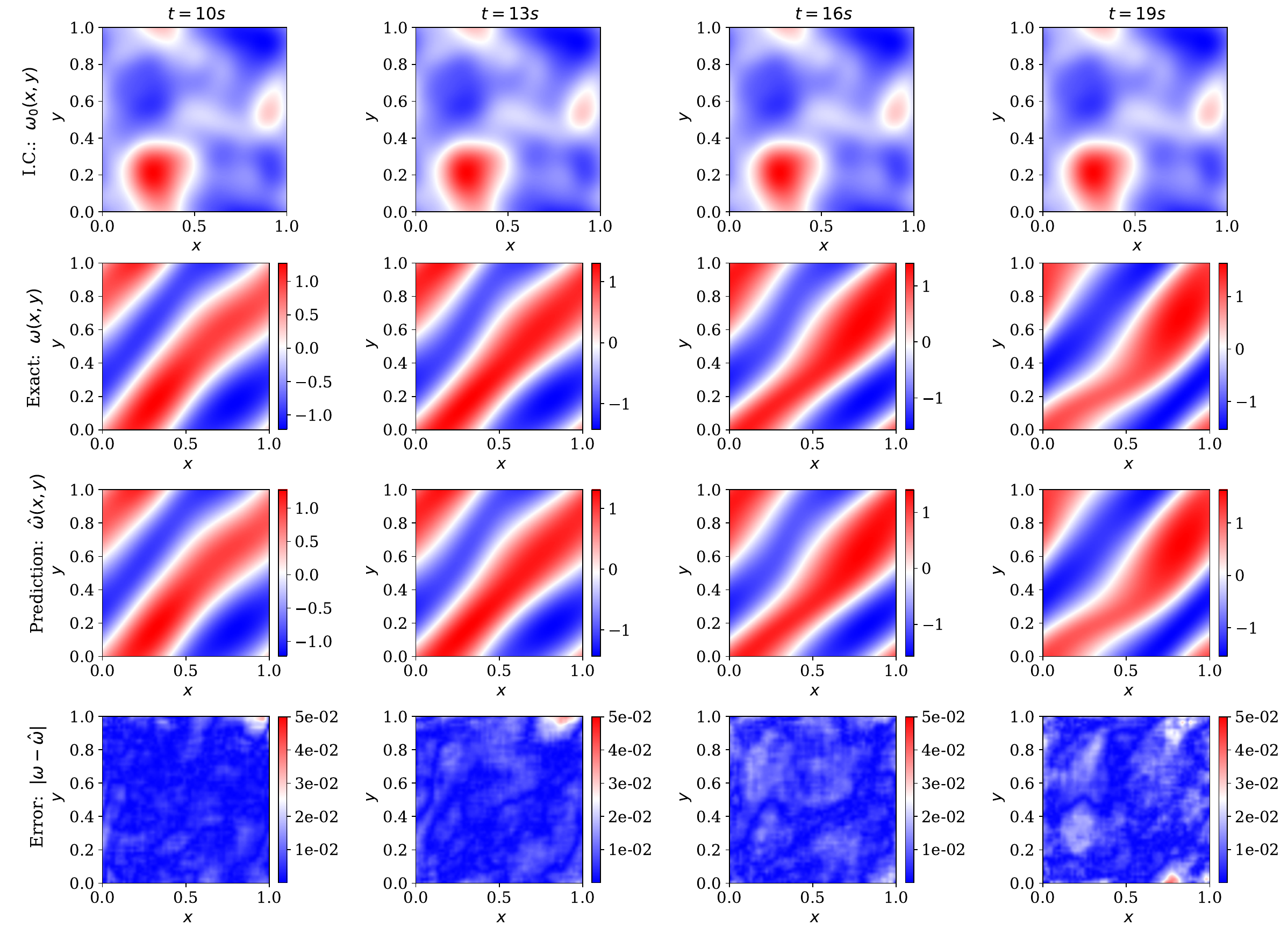}
	\caption{\; Navier-Stokes equations with a spatial resolution of $64 \times 64$. Vorticity data at several initial time steps $t \in [0, 10]$ are used to predict vorticity data and error maps at subsequent time steps $t \in (11, 20]$.}
	\label{fig13}
\end{figure}

\textbf{Results}: The prediction results for vorticity field at $t \in (10, 20]$ are given in Figure \ref{fig13}, along with the mean $L_2$ relative test error in Table \ref{tab3}. The results show that the prediction error of the proposed U-WNO is second only to MWT, lower than WNO and FNO, and less than $1/3$ of the error of WNO. It can be speculated that the proposed U-WNO can learn the operators of the 2D time-dependent PDE with higher accuracy than WNO. Table \ref{tab5} shows the long-term predictions of the Navier-Stokes equations by WNO and U-WNO. Although the training data only covers $t \in (0, 10]$, U-WNO still maintains a low error in the long-term prediction of $t \in (10, 40]$, indicating that it has some extrapolation ability.

\begin{table}[ht!]  
\centering  
	\caption{\; The mean $L_2$ relative error of the Navier-Stokes equations for long time predictions} 
	\label{tab5}  
\begin{tabular}{cccccc}  
\hline  
Network architecture &$t \in (10, 20]$ & $t \in (10, 25]$ & $t \in (10, 30]$ & $t \in (10, 35]$ & $t \in (10, 40]$ \\ \hline  
U-WNO               & 0.96\%          & 1.40\%          & 2.77\%           & 5.21\%          & 8.06\%           \\
WNO                & 3.16\%           & 4.34\%          & 7.24\%          & 11.60\%          & 16.28\%         \\ \hline  
\end{tabular}  
\end{table} 

\section{Conclusions}

In this paper, we propose the U-WNO neural operator, which integrates a U-Net path and an improved wavelet layer with an adaptive activation function. Unlike the baseline wavelet neural operator (WNO) that parameterizes only the integral kernel on coarse-grained discrete wavelet transform (DWT) scales, we introduce a U-Net pathway into the wavelet layer to enhance high-frequency feature extraction. By incorporating trainable slope parameters in the activation function, U-WNO mitigates spectral bias and further improves high-frequency representation capabilities. Additionally, we design residual connections inspired by neural network skip pathways to ensure spatial information completeness in each enhanced wavelet layer. Extensive numerical experiments compare U-WNO with WNO and other state-of-the-art operators (DeepONet, FNO, MWT). Results demonstrate that U-WNO achieves the lowest $L_2$ relative errors across multiple benchmarks and effectively learns nonlinear operators in complex geometries. Future work will focus on applying U-WNO to real-world engineering problems and developing a hybrid training framework that combines data-driven methods with physics-informed constraints.

%In this paper, we propose U-WNO, a neural operator based on U-Net, and an improved wavelet layer with the Adaptive Activation Function. The WNO parameterizes only the integral kernel on a higher scale of DWT. Therefore, we add a U-Net path in the wavelet layer to enhance the model's performance in capturing high-frequency information. By introducing the Adaptive Activation Function into U-WNO, we overcome the Spectral Bias and further enhance U-WNO's ability to express high-frequency information. Furthermore, we design a residual shortcut based on the residual connection of NN to ensure that each improved wavelet layer of U-WNO receives complete location information, thus enhancing U-WNO's expressiveness in the spatial domain. In the numerical experiments, we compare U-WNO with WNO, as well as other popular neural operators like DeepONet, FNO, and MWT. The proposed U-WNO demonstrates excellent performance, achieving the lowest error on multiple examples. U-WNO also proves capable of learning nonlinear operators in complex geometric domains. In the future, we plan to apply U-WNO to real engineering problems to assess its advantages and disadvantages. Additionally, since U-WNO, similar to WNO, is data-driven, its training still relies on data obtained from traditional solvers. To further improve U-WNO, we intend to adopt a physical information-driven approach, enabling it to learn from both data-driven methods and control physics.

\textbf{Acknowledgements}. \emph{The authors are very much indebted to the referees for their constructive comments and valuable suggestions}.

%\textbf{Conflict of Interest Statement}\\
%\indent This study does not have any conflicts to disclose.

\textbf{Code Availability Statement}\\
\indent The code used to generate the findings of this study is openly available in GitHub at \url{https://doi.org/10.5281/zenodo.15612292}.

%\textbf{ORCID}\\
%\indent Hou-biao Li \url{https://orcid.org/0000-0002-7268-307X.}

\vspace{-6mm}
\section*{}
\bibliographystyle{unsrt}
\bibliography{ref.bib}
%\end{thebibliography}
\end{document}